\documentclass[aap,MSNbibl,citesort,dvips]{arximspdf}
\usepackage{mathrsfs}

\doi{10.1214/11-AAP812} 
\volume{22}
\issue{5}
\pubyear{2012}
\firstpage{1745}
\lastpage{1777}

\makeatletter

\newcommand{\bin}{\operatorname{Bin}}

\newcommand{\eqdist}{\stackrel{d}{=}}
\newcommand{\one}{\mathbf{1}}

\newtheorem{theorem}{Theorem}[section]
\newtheorem{prop}{Proposition}[section]
\newtheorem{Lemma}{Lemma}[section]
\newtheorem{cor}{Corollary}[section]

\newproclaim{Example}{Example}
\newproclaim{Plan}{Plan of the paper}
\newproclaim{algorithmic}{An algorithmic point of view}
\newproclaim{Background}{Background and previous work}
\newproclaim{Discussion}{Discussion and remarks about the assumptions}
\newproclaim{Remarks}{Remarks}

\makeatother

\begin{document}
\begin{frontmatter}

\title{The total path length of split trees}
\runtitle{Path length of split trees}

\begin{aug}
\author[A]{\fnms{Nicolas} \snm{Broutin}\corref{}\ead[label=e1]{nicolas.broutin@inria.fr}}
\and
\author[B]{\fnms{Cecilia} \snm{Holmgren}\ead[label=e2]{cecilia@math.uu.se}}
\runauthor{N. Broutin and C. Holmgren}
\affiliation{Inria, and Inria and Cambridge University}
\address[A]{Projet Algorithms\\
Inria Paris-Rocquencourt \\
Domaine de Voluceau\\
78153 Le Chesnay\\
France\\
\printead{e1}}
\address[B]{DPMMS\\
Centre for Mathematical Sciences\\
Wilberforce Road\\
Cambridge CB3 0WA\\
United Kingdom\\
\printead{e2}} 
\end{aug}

\received{\smonth{2} \syear{2011}}
\revised{\smonth{9} \syear{2011}}

%
\begin{abstract}
We consider the model of random trees introduced by Devroye
[\textit{SIAM J. Comput.} \textbf{28} (1999) 409--432]. The model
encompasses many important randomized algorithms and data structures.
The pieces of data (items) are stored in a randomized fashion in the
nodes of a~tree. The total path length (sum of depths of the items) is
a natural measure of the efficiency of the algorithm/data structure.
Using renewal theory, we prove convergence in distribution of the total
path length toward a distribution characterized uniquely by a fixed
point equation. Our result covers, using a unified approach, many data
structures such as binary search trees, $m$-ary search trees, quad
trees, median-of-$(2k+1)$ trees, and simplex trees.
%
\end{abstract}

%

\end{frontmatter}

\section{Introduction}\label{intro}

In this paper we investigate the total path length, that is, sum of all
depths, of random split trees defined by Devroye~\cite{Devroye1998}
(we will
be more precise shortly).
Split trees model a large class of efficient data structures or sorting
algorithms.
Some important examples of split trees are
binary search trees (which are also the representation of Quicksort)
\cite{Hoare1962}, $m$-ary search trees~\cite{Pyke1965},
quad trees~\cite{FiBe1974}, median-of-$(2k+1)$ trees~\cite{Bell1965},
simplex trees; all these are covered by the results in this document.
The case of tries~\cite{Fredkin1960} and digital search trees
\cite{CoEv1970} is also important in practice~\cite{Szpankowski2001};
however, their treatment necessitates different tools, and we leave
this case for later studies.

The magnitude of the depths in tree data structures naturally
influences their efficiency; in the case where the tree represents the
branching choices made by an algorithm, the depths are related to the
running time of the algorithm. In this sense, the sum of the depths is
a natural and important measure of the efficiency of tree data
structures or sorting algorithms. 



The path length of tree data structures has been studied by many
authors, but in most cases the analyses and proofs are very much tied
to a specific case. The main result of this study is to prove that for
a large class of split trees, the total path length converges in
distribution to a random variable characterized by some fixed point
equation. In that sense our result extends the earlier studies of
R{\"o}sler~\cite{Roesler1991,Roesler2001} and Neininger and
R{\"u}schendorf~\cite{NeRu1999} who used the so-called contraction
method to show convergence in distribution of the total path length for
the specific examples of the binary search trees, the
median-of-$(2k+1)$ trees and quad trees. Our method actually relies on
previous work of Neininger and R{\"u}schendorf~\cite{NeRu1999} who gave
a limit theorem for the path
length of general split trees, under the assumption that the mean
satisfies some precise asymptotic form, which we prove.
\begin{Plan*}
In Section~\ref{secbackground},
we introduce the model of split trees of Devroye~\cite{Devroye1998}.
We also
discuss previous work on the path length and similar topics. This is
also the place where we state our main result, Theorem~\ref{main}.

In Section~\ref{seccontraction}, we explain our general approach,
which relies heavily on previous work by Neininger and R{\"u}schendorf
\cite{NeRu1999}. These authors stated a~general condition for
convergence in distribution of the path length, and our contribution is
to prove that it indeed holds for a large class of split trees. So
Section~\ref{seccontraction} is included so that the reader has a
general view of the argument.

Once we have stated the precise condition in Section \ref
{seccontraction}, we will move on to explaining our approach to
proving it in Section~\ref{secaverage}. Finally, in Section
\ref{secconcl} we discuss extensions of our results.
\end{Plan*}

\section{Split trees and path length: Notation and background}
\label{secbackground}

We introduce the split tree model of Devroye~\cite{Devroye1998}.
Consider an
infinite rooted $b$-ary tree (every node has $b$ children). The nodes
are identified with the set of finite words on an alphabet with $b$
letters, $\mathcal U=\bigcup_{n\ge0}\{1,\ldots, b\}^n$. The root is
represented by the empty word~$\varnothing$. We write $u\preceq v$ to
denote that $u$ is an ancestor of $v$ (as words, $u$ is a prefix of
$v$). In particular, for the empty word $\varnothing$, we have
$\varnothing\preceq v$ for any $v\in\mathcal U$.

A split tree $T^n$ of \textit{cardinality} $n$ is constructed by
distributing $n$ items (pieces of data) to the nodes $u\in\mathcal U$.
To describe the tree, it suffices to define the number of items $n_u$
in the subtree rooted at any node $u\in\mathcal U$. The tree~$T^n$ is
then defined as the smallest relevant tree, that is, the subset of
nodes~$u$ such that $n_u>0$ (which is indeed a tree).

In the model, internal nodes all contain $s_0\ge0$ items, and external
nodes can contain up to $s$ items.
The construction then resembles a divide-and-conquer procedure, where
the partitioning pattern depends on a random vector of proportions. Let
$\mathcal V=(V_1,\ldots,V_b)$ satisfy $V_i\ge0$ and $\sum_i V_i=1$;
each node $u\in\mathcal U$ receives an independent copy $\mathcal V_u$
of the random vector $\mathcal V$. In the following, we always assume
that ${\mathbf P}(\exists i\dvtx V_i=1)<1$.
We can now describe $(n_u, u\in\mathcal U)$.
The tree contains $n$ items, and we naturally have $n_\varnothing$. The
split procedure is\vadjust{\goodbreak} then carried on from parent to children as long as
$n_v>s$. Given the cardinality $n_v$ and the split vector $\mathcal
{V}_v=(V_1,V_2,\ldots,V_b)$ of $v$, the cardinalities
$(n_{v_1},n_{v_2},\ldots,n_{v_b})$ of the $b$ subtrees rooted at
$v_1,v_2,\ldots,v_b$ are distributed as
%
\begin{equation}\label{mult}
\operatorname{Mult}(n_v-s_0-bs_1,V_1,V_2,\ldots, V_b)+(s_1,s_1,\ldots, s_1),
\end{equation}
where $0\leq s$ and $0\leq bs_1\leq s+1-s_0$.



Depending on the choice of parameters $s_0, s_1, s$ and the
distribution of $\mathcal V=(V_1, \ldots, V_b)$, many important data
structures may be modeled, such as binary search trees, $m$-ary search
trees, median-of-$(2k+1)$ trees, quad trees, simplex trees
(see~\cite{Devroye1998}). To make sure that the model is clear and to give
a hint of the wide applicability of the model, we illustrate the
construction with two canonical examples.
%
\begin{Example}[(Binary search tree)]\label{Ex1} 
The binary search tree is one of the most common data structures for
sorted data. Here we assume that the data set is $\{1,\ldots,n\}$. A
first (uniformly) random key is drawn $\sigma_1$, and stored at the
root of a binary tree. The remaining keys are then divided into two
subgroups, depending on whether they are smaller or larger than $\sigma
_1$. The left and right subtrees are then binary search trees built
from the two subgroups $\{i\dvtx i<\sigma_1\}$ and $\{i\dvtx i>\sigma_1\}$,
respectively.
The sizes of the two subtrees of the root are $\sigma_1-1$ and
$n-\sigma
_1$. One easily verifies that, since $\sigma_1$ is uniform in $\{
1,2,\ldots,n\}$, one has
\[
(\sigma_1-1, n-\sigma_1)\eqdist\operatorname{Mult}(n-1; U, 1-U),
\]
where $U$ is a uniform $U(0,1)$ random variable.
Thus, a binary search tree can be described as a split tree with
parameters $b=2$, $s_0=1$, $s=1$, $s_1=0$ and $\mathcal V$ is
distributed as $(U,1-U)$ for $U$ a random variable uniform on $[0,1]$.
\end{Example}

\begin{Example}[(Digital trees or tries)]\label{Ex2}
We are given $n$ (infinite) strings $X_1,\break\ldots, X_n$ on the alphabet
$\{1,\ldots,b\}$. The strings are drawn independently, and the symbols
of every string are also independent with distribution on
$\{1,\ldots,b\}$ given by $p_1,\ldots, p_b$. Each string naturally
corresponds to an infinite path in the infinite complete $b$-ary tree,
where the sequence of symbols indicates the sequence of directions to
take as one walks away from the root. The trie is then defined as the
smallest tree so that all the paths corresponding to the infinite
strings are eventually distinguished; that is, for every string~$X_i$,
there exists a node $u$ in the tree such that $X_i$ is the only string
with $u\preceq X_i$. The internal nodes store no data; each leaf stores
a unique string. In this case, $n_v$ is the number of strings that have
prefix $v$, and one clearly has for the children of the root
\[
(n_{1},\ldots, n_b)\eqdist\operatorname{Mult}(n; p_1,\ldots, p_b).
\]
The trie is thus a random split tree with parameters $s=1, s_0=s_1=0$
and $\mathcal V=(p_1,p_2,\ldots, p_b)$ almost surely.\vadjust{\goodbreak}
\end{Example}
%
%
\begin{algorithmic*}
Rather than using the
divide-and-conquer description above, the random trees may be
equivalently defined using incremental insertion of data items into an
initially empty data structure. The items are labeled using $\{
1,2,\ldots, n\}$ in the order of insertion. Initially, $n_u=0$ for
every $u\in U$. We first sample the i.i.d. copies of $\mathcal V$ that
are assigned to the nodes $u\in\mathcal U$.
%
\begin{itemize}
\item Upon insertion, an item first trickles down along a random path
from the root until it finds a leaf (i.e., a node $u$ such that all its
children $u_1,\ldots, u_b$ satisfy $n_{u_i}=0$). If the path currently
corresponds to a word $v\in\mathcal U$, and $v$ is not a leaf, then it
is extended to $v_i$, the $i$th child of $v$ with probability~$V_i$,
where $(V_1,\ldots, V_b)$ is the copy of $\mathcal V$ associated with $v$.
\item When the first phase is finished, the item is stored in a leaf,
say $v$. The leaves can contain up to $s$ items. So if $n_v<s$ (before
the insertion), then the item is stored at $v$, and all the $n_u$ for
$u\preceq v$ are updated.
\item If $n_v=s$, there is no space for the new item at $v$. With the
new item, we formally have $n_v=s+1$. In this case, $s_0$ of these
$s+1$ items are randomly chosen to remain at $v$ while the other
$s+1-s_0$ are distributed among the children $v_1,\ldots,v_b$ of $v$.
Each child receives $s_1$ items chosen at random. The remaining
$s+1-s_0-b s_1$ each choose (independently) a child $v_i$ at random
with probability $V_i$, where $(V_1,\ldots, V_b)$ is the copy of
$\mathcal V$ at node $v$. If $s_1=s_0=0$, it may happen that all $s+1$
items now lie at one child $v_i$, in which case the scheme is repeated
until a stable position is found. [This happens with probability 1,
since ${\mathbf P}(\exists i\dvtx V_i=1)<1$.] This last step is the reason
why an
item may move down when a further item is inserted.
\end{itemize}

The properties of the multinomial distribution ensure that the tree
$T^n$ obtained in this way has the correct distribution (see
\cite{Devroye1998} for details).

In the present case we can assume without loss of generality that the
components of $\mathcal{V}$ are identically distributed; applying a
random permutation to the components would leave the path length
unchanged. We now let $V$ denote a uniformly random component of
$\mathcal V$. So, for instance, $\mathbf{E}[V]=1/b$ and ${\mathbf
P}(V=1)<1/b$ by our
assumption that ${\mathbf P}(\exists i\dvtx V_i=1)<1$.
\end{algorithmic*}
\begin{Background*}
The labeling of the items induced by the algorithm above is interesting
for the analysis. Let $D_i$ be the depth of the item labeled $i$ when
all $n$ items have been inserted. Then, the total path length is
\[
\Psi(T^n)=\sum_{i=1}^n D_i.
\]
The analysis of the depth $D_n$ of the last item $n$ is thus tightly
related to the analysis of $\Psi(T_n)$, and yet is much simpler since
it avoids the intricate dependence between the $D_i$.\vadjust{\goodbreak}
Devroye~\cite{Devroye1998} proved a weak law of large numbers and a central
limit theorem for $D_n$ in general split trees.
Let $\Delta$ be a~component of $(V_1,\ldots,V_b)$ picked with
probability proportional to its size; that is, given $(V_1,\ldots,V_b)$,
let $\Delta=V_j$ with probability $V_j$. We write
%
\begin{eqnarray}\label{mu}
\mu&:=&\mathbf{E}[-\ln\Delta]=b\mathbf{E}[-V\ln{V}]
\quad\mbox{and}\nonumber\\[-8pt]\\[-8pt]
\sigma^2&:=&\operatorname{\mathbf{Var}}(\ln\Delta)=b\mathbf
{E}[V\ln^2 V]-\mu^2.\nonumber
\end{eqnarray}
Note that $\mu\in(0,\infty)$ and $\sigma<\infty$. Then $D_n/\ln n$
converges in probability to~$\mu^{-1}$, and $\mathbf{E}[D_n]/\ln
n\rightarrow
\mu^{-1}$ (Devroye assumed that ${\mathbf P}(V=1)=0$, but this
assumption can be
relaxed as long as $V$ satisfies ${\mathbf P}(V=1)<1/b$; this is done using
trees in which edges are weighted by geometric random variables
(see, e.g.,~\cite{BrDe2006,BrDeMc2008})). If we also have $\sigma>0$, then
\[
\frac{D_n-\mu^{-1} \ln n}{\sqrt{\sigma^2\mu^{-3}\ln n}}\to
\mathcal N(0,1)
\]
in distribution where $\mathcal N(0,1)$ denotes the standard Normal
distribution. Note that $\sigma>0$ precisely when $V$ is not
monoatomic, that is, if $bV\ne1$ with positive probability.

The total path length $\Psi(T^n)$ itself has been extensively studied
for specific cases of split trees. The first moment follows from that
of $D_n$ since
\[
\mathbf{E}[\Psi(T^n)]=\sum_{i=1}^n \mathbf{E}[D_i].
\]
For instance, in the binary search tree, we have~\cite{Hibbard1962}
%
\begin{equation}\label{bin}
\mathbf{E}[\Psi^{\mathrm{BST}}(T^n)]=2n\ln{n}+n(2\gamma-4)+2\ln
{n}+2\gamma+1+
\mathcal O (n^{-1}),
\end{equation}
where $\gamma$ is Euler's constant.
For higher moments and the distribution of~$\Psi(T^n)$, one needs to
carefully take the dependence in the terms of the sum into account.
Most studies of this type concern the model of binary search tree, or
equivalently the cost of quicksort (e.g.,
\cite{Regnier1989,Roesler1991,FiJa2002,FiJa2002,TaHa1995}).
Let
%
\begin{equation}
Y_n:=\frac{\Psi^{\mathrm{BST}} (T^n)-\mathbf
{E}[\Psi^{\mathrm{BST}}(T^n)]}{n}.
\end{equation}
Using martingale arguments, R{\'e}gnier~\cite{Regnier1989} showed that $Y_n$
converges in distribution to a random variable $Y$.
R{\"o}sler~\cite{Roesler1991}
showed that $Y$ is satisfying the following distributional equality:
%
\begin{equation}\label{regnier}
Y\eqdist UY+(1-U)Y^*+C(U),
\end{equation}
where $C(u):=2u\ln u+2(1-u)\ln(1-u) +1$, $U$ is uniform on $[0,1]$, $Y$
and $Y^*\eqdist Y$ are independent. He also proved that the stochastic
equality in (\ref{regnier}) actually characterizes the distribution of
$Y$: there exists a unique solution $Y$ of (\ref{regnier}) such that
$\mathbf{E}[Y]=0$ and\vadjust{\goodbreak} $\operatorname{\mathbf{Var}}(Y)<\infty$. The
distribution of $Y$ is usually
called the quicksort distribution. Properties of $Y$ and the rate of
convergence of $Y_n$ to $Y$ are studied in
\cite{Roesler1991,FiJa2001,FiJa2002,TaHa1995}.

The aim of the present study is to prove that the path length exhibits
a~similar asymptotic behavior regardless of the precise model of split
tree: 

\begin{theorem}\label{main}Let $ \Psi(T^n) $ be the total path length
in a
general split tree with split vector $\mathcal{V}=(V_1,\ldots,V_b)$.
Suppose that ${\mathbf P}(\exists i\dvtx V_i=1)<1$. Let
\[
X_n:=\frac{\Psi(T^n)-\mathbf{E}[\Psi(T^n)]}n \quad\mbox{and}\quad
C(\mathcal{V})=1+\frac1 \mu\sum_{i=1}^{b}V_i\ln V_i.
\]
If $C(\mathcal V)\ne0$ with positive probability, then $X_n\to X$ in
distribution, where $X$ is the unique solution of the fixed point equation
\[
X\eqdist\sum_{k=1}^{b}V_kX^{(k)}+C(\mathcal{V}),
\]
satisfying $\mathbf{E}[X]=0$ and $\operatorname{\mathbf
{Var}}(X)<\infty$.
Furthermore, exponential moments of $X_n$ exist and converge $\mathbf
{E}[e^{\lambda X_n}]\rightarrow\mathbf{E}[e^{\lambda X}]$ for any
$\lambda\in
\mathbb R$.
\end{theorem}

As mentioned in the \hyperref[intro]{Introduction}, Neininger and R{\"u}schendorf
\cite{NeRu1999} proved a version of Theorem~\ref{main} conditional on
the type of asymptotic expansion for ${\mathbf E}\Psi(T^n)$; our
contribution is
to prove that this expansion indeed holds (Theorem~\ref{expectedpath}),
which implies the unconditional version stated in Theorem~\ref{main}.

We have recently been informed that, based on a Markov chain
representation of Bruhn~\cite{Bruhn1996a} and coupling arguments,
Munsonius~\cite{Munsonius2011a} has shown a~result similar to our
Theorem~\ref{main}
in the special case when the distribution of $V$ has a density with
respect to Lebesgue measure.
\end{Background*}
\begin{Discussion*}
(i) When the split vector $\mathcal V$ is deterministic, that
is, $\mathcal V$ is a permutation of some fixed vector $(p_1,\ldots,
p_b)$, the cost function $C(\mathcal V)=0$. Such a split tree is a
digital tree~\cite{Szpankowski2001}. In some sense, part of
Theorem~\ref{main} still holds, but the limit $X$ is trivial since
$X=0$ almost surely. The renormalization is actually too strong, since
the variance in this case should be of order $n\log n$, rather than
$n^2$ [and order $n$ in the special case when $b\mathcal
V=(1,\ldots,1)$]. The total path length for binary tries has been
treated by Jacquet and R\'egnier~\cite{JaRe1988}. They showed that the
variance of $\Psi(T_n)$ is of order $\mathcal{O}(n)$ if $p=q$ and of
order $\mathcal{O}(n\log n)$ if $p \neq q$ and that the path length is
asymptotically normal. Schachinger~\cite{Schachinger2004} showed that,
for tries
with a general branch factor, the variance of the total path length for
general tries is $\mathcal
O(n\log^2 n)$. See also~\cite{KiPrSz1994a,KiPrSz1989a}.%

\vspace*{-6pt}

\begin{longlist}[(iii)]
\item[(ii)] In general, in the case of digital trees [when
$C(\mathcal V)=0$], it is expected that under the correct rescaling the
limit distribution should be normal. Neininger and R{\"u}schendorf
\cite{NeRu2004} gave a general
conditions under which limit distributions\vadjust{\goodbreak} are Gaussian. The case of
the binary tries is one example when this theorem can be applied as an
alternative proof to the method in~\cite{JaRe1988}. In general, to
apply the result in~\cite{NeRu2004} one needs to have approximations
for the first two moments of the path length. This is the reason why we
report the analysis of this case: a lot more work is required to
estimate the variance to the correct order.

\item[(iii)] It might seem at first that one should have $C(\mathcal
V)=0$ when $\ln V$ is lattice (trie case). However, one can easily
construct examples with $C(\mathcal V)\ne0$ and $\ln V$ lattice: for
instance, take $b=5$ and $\mathcal V$ a random permutation of either
$(1/2,1/8,1/8,1/8, 1/8)$ or $(1/2,1/4,1/4,0,0)$, each with
probability $1/2$.

\item[(iv)] Note, although it might come as a surprise since our
main tool is renewal theory, Theorem~\ref{main} does not require any
condition on arithmetic properties related to the vector $(V_1,\ldots,
V_b)$. In particular, it holds whether $-\ln V$ is lattice or not.
However, the behavior of the average path length does depend on
arithmetic properties of $\ln V$; see Theorem~\ref{expectedpath} later
for details.

\item[(v)] Note that the limit fixed equation only depends on
$\mathcal V$, so in particular, the limit distribution $X$ does not
depend on the parameters $s,s_0$ or $s_1$. However, the average
$\mathbf{E}[\Psi(T^n)]$ should clearly depend on these parameters,
although we do
not prove it formally.

\item[(vi)] For the sake of simplicity, we cover only trees with
bounded degree, which is usually the case for trees representing data
structures. The path length of recursive trees, which do not have
bounded degree, has been studied by~\cite{Mahmoud1991,DoFi1999}.
\end{longlist}
\end{Discussion*}


\section{The contraction method for path length}\label{seccontraction}

The condition stated by Neininger and R{\"u}schendorf~\cite{NeRu1999}
to ensure weak convergence of the path length concerns the asymptotics
of the average path length. More precisely, if one has, for some
constant $\varsigma$,
%
\begin{equation}\label{eqexpectedcond}
\mathbf{E}[\Psi(T^n)]=\mu^{-1} n\ln n +\varsigma n +o(n)
\end{equation}
and ${\mathbf P}(C(\mathcal V)\ne0)>0$, then Theorem 5.1 of~\cite{NeRu1999}
ensures that $X_n\to X$ in distribution.
The purpose of this section is to explain why these conditions are
sufficient to prove Theorem~\ref{main}. In particular, we give the
necessary background about the contraction method, and we explain the
general approach that has been devised in~\cite{NeRu1999}. This section
is included only to put our result in context, and no new result is
proved with respect to the contraction method.

Note first that (\ref{eqexpectedcond}) holds in the case of binary
search trees (\ref{bin}). Recall that~$D_i$ is the depth of the $i$th
item in the construction where items are inserted one after another. It
is not difficult to deduce from the results on~$D_i$ by Devroye
\cite{Devroye1998} that
\[
\mathbf{E}[\Psi(T^n)]=\mu^{-1} n \ln n + n q(n)
\]
with $q(n)=o(\ln n)$ (see Theorem 2.3 of~\cite{Holmgren2009} for a
formal proof). So proving~(\ref{eqexpectedcond}) reduces to proving
that $q(n)\to\varsigma$ as $n\to\infty$. Our contribution is to prove
that this is indeed the case as soon as the random variable $V$ is such
that $-\ln V$ is not lattice, that is, there is no $a\in\mathbb R$
such that $-\ln V\in a\mathbb Z$ almost surely. In the following, we let
\[
d:=\sup\{a\ge0\dvtx {\mathbf P}(\ln V\in a\mathbb Z)=1\},
\]
so that $d$ is the span of the lattice when $d>0$ and $\ln V$ is
nonlattice when $d=0$.
More precisely, we prove:
%
\begin{theorem}\label{expectedpath}
The expected value of the total path length $\Psi(T^n) $ exhibits the
following asymptotics, as $n\to\infty$:
%
\begin{equation}\label{exptotal1h}
\mathbf{E}[\Psi{(T^{n})}]=\mu^{-1} n\ln n+n \varpi(\ln n)+o(n),
\end{equation}
where $\mu$ is the constant in (\ref{mu}) and $\varpi$ is a continuous
periodic function of period~$d$. In particular, if $\ln V$ is not
lattice, then $d=0$ and $\varpi$ is constant.
\end{theorem}

If $\ln V$ is nonlattice, then Theorem~\ref{expectedpath} and
Theorem 5.1 of~\cite{NeRu1999} together prove Theorem~\ref{main}. If
the random variable $\ln V$ is lattice with span $d$, then
Theorem~\ref{expectedpath} implies that $q(n)=\varpi(\ln n)+o(1)$ as
$n\to\infty$, where $\varpi$ is $d$-periodic. So it seems that
Theorem~\ref{expectedpath} does not permit to conclude along the
arguments by Neininger and R{\"u}schendorf~\cite{NeRu1999}. However,
the techniques in~\cite{NeRu1999} only require convergence of the
coefficients of a certain recursive equation; this fact was used in
\cite{NeRu2004} to deal with certain cases involving oscillations.




We now move on to the approach developed by Neininger and
R{\"u}schendorf \mbox{\cite{NeRu1999,NeRu2004}}. Let $\overline
n=(n_1,\ldots,n_b)$ denote the vector of cardinalities of the children
of the root. Then we have, for $n>s$,
\[
\Psi(T^n)\stackrel{d}=\sum_{i=1}^{b}\Psi_i (T^{n_i})+n-s_0,
\]
where $\Psi_i (T^{n_i})$ are copies of $\Psi(T^{n_i})$ that are
independent conditional on $(n_1,\ldots, n_b)$. Introducing the
normalized total path length
%
\begin{equation}\label{total}
X_n:=\frac{\Psi(T^n)-\mathbf{E}[\Psi(T^n)]}{n},
\end{equation}
we can rewrite the distributional identity above as
\[
X_n:=\sum_{i=1}^{b}\frac{n_i}{n}X_{n_i}+C_n(\overline{n}),
\]
where
\[
C_n(\overline{n}):=1-\frac{s_0}{n}-\frac{\mathbf{E}[\Psi
(T^n)]}n+\sum
_{i=1}^b\frac{\mathbf{E}[\Psi(T^{n_i})]}n\vadjust{\goodbreak}
\]
and $X_{n_i}$, $i\in\{1,\ldots,b\}$, are independent conditional on
$(n_1,\ldots,n_b)$.
By definition, the vector of cardinalities $\overline{n}$ is
$\operatorname{Mult}(n-s_0-bs_1,V_1,V_2,\ldots, V_b)+(s_1,s_1,\ldots,
s_1)$ so that
%
\begin{equation}\label{eqlimitsplit}
\biggl(\frac{n_1} n,\frac{n_2}n, \ldots, \frac{n_b}n\biggr)\rightarrow
\mathcal
V_\sigma=(V_1,V_2,\ldots, V_b),
\end{equation}
almost surely as $n\to\infty$. This is where (\ref{eqexpectedcond})
comes into play: it ensures that the cost $C_n(\overline n)$ (the
``toll function'') in the recursive distributional equation does
converge (in distribution) as $n\to\infty$. Indeed
\begin{eqnarray*}
C_n(\overline n)
&=&1+\frac1 n \sum_{i=1}^b \mathbf{E}[\Psi(T^{n_i})]-\frac{\mathbf
{E}[\Psi(T^n)]}n -
\frac{s_0}n\\
&=&1+\frac1 \mu\sum_{i=1}^b \frac{n_i}n \ln\frac{n_i}n + \frac1
\mu
\Biggl(\sum_{i=1}^b \frac{n_i}n \varpi(\ln n_i) - \varpi(\ln n)\Biggr)+o(1).
\end{eqnarray*}
Now, by (\ref{eqlimitsplit}) and the continuity of $\varpi$, it
follows that
%
\begin{eqnarray}\label{eqtoll}\qquad
C_n(\overline n)
&=&1+\frac1 \mu\sum_{i=1}^b \frac{n_i}n \ln\frac{n_i}n + \frac1
\mu
\Biggl(\sum_{i=1}^b \frac{n_i}n \varpi(\ln n + \ln V_i) - \varpi(\ln
n)\Biggr)+o(1)\nonumber\\[-8pt]\\[-8pt]
&=&1+\frac1 \mu\sum_{i=1}^b V_i \ln V_i +o(1),\nonumber
\end{eqnarray}
since $\varpi$ is $d$-periodic and $\ln V_i\in d \mathbb Z$ by
assumption (if $d=0$, $\varphi$ is constant and the claim also holds).
Note that, apart from (\ref{eqlimitsplit}), only asymptotics for the
first moments are required for (\ref{eqtoll}) to hold. Together
(\ref{eqlimitsplit}) and (\ref{eqtoll}) suggest that if $X_n$
converges in
distribution to some limit $X$, then
$X$ should satisfy the following fixed point equation:
%
\begin{equation}\label{fixsplit}
X\eqdist\sum_{k=1}^{b}V_kX^{(k)}+C(\mathcal{V})
\qquad\mbox{where } C(\mathcal V)=1+\frac1 \mu\sum_{i=1}^b V_i \ln V_i,
\end{equation}
and $X^{(k)}$ are independent and identically distributed copies of $X$.

The point of the contraction method is to make the previous arguments
rigorous, that is, to show that if the coefficients $C_n(\overline n)$
do converge, then~(\ref{fixsplit}) has a unique solution $X$ and that
$X_n\to X$ in distribution; this is precisely what was done in
\mbox{\cite{NeRu1999,NeRu2004}}. This is done by proving that the recursive map
defined by (\ref{fixsplit}) is a contraction in a suitable space of
probability measures~\cite{Roesler1991,Roesler1992,RaRu1995}. We now
expose the lines of the arguments to show the extent of the results
that follow from the mere convergence of the coefficients~$C_n(\overline n)$.
(We claim no novelty.)

Let $\mathscr M_2$ be the set of probability measures with a finite
second moment. For a random variable $X$, we write $\mathcal D(X)$ for
its law. For $\phi\in\mathscr M_2$ and $X$ a random variable with law
$\mathcal D(X)=\phi$, define the $L^2$-norm by $\| X\|_2=\mathbf
{E}[X^2]^{1/2}$. We can then define a metric $d_2$ on $\mathscr M_2$ (the
Mallow metric): for $\phi,\varphi\in\mathscr M_2$, let
%
\begin{equation}
d_2(\phi,\varphi):={\inf}\| X-Y \|_2,
\end{equation}
where the range of the infimum is the set of couples $(X,Y)$ with
marginal distributions $\mathcal D(X)=\phi$ and $\mathcal D(Y)=\varphi$.
For simplicity we write $d_2(X,Y)=d_2(\phi,\varphi)$ for random
variables $X$ and $Y$, but note that this only depends on the marginal
distributions $\phi$ and $\varphi$.
Convergence of $\phi_n$ to $\phi$ in $(\mathscr M_2,d_2)$ is equivalent
to weak convergence with convergence of the second moment~\cite{RaRu1995}:
%
\begin{equation}\label{Mallow}
\phi_n\stackrel w\rightarrow\phi \quad\mbox{and}\quad \int x^2 \,d\phi
_n(x)\rightarrow\int x^2 \,d\phi(x).
\end{equation}

Let $\mathscr M_2^0$ be the subset of $\mathscr M_2$ containing
distributions $\phi$ such that $\int\! x \,d\phi(x)=0$. Define the operator
$T\dvtx \mathscr M_2^0\to\mathscr M_2^0$. For a distribution $\phi\in
\mathscr M_2^0$, let $T(\phi)$ be the distribution of the random
variable given by
\[
\sum_{1\le k\le b}V_kZ^{(k)}+C(\mathcal{V}),
\]
where $Z^{(i)}$ are i.i.d. random variables with distribution $\phi$.
Then, calculations similar to that in the proof of Lemma 3.2 in
\cite{NeRu1999} yield
\begin{eqnarray*}
d_2(T(X),T(Y))
&\leq&\sum_{1\le i\le b}\mathbf{E}[V_i^2]\cdot d_2(X,Y)\\
&=& b \mathbf{E}[V^2] \cdot d_2(X,Y).
\end{eqnarray*}
Since $b\mathbf{E}[V^2]<1$ the operator $T$ is a contraction in
$(\mathscr
M_2^0, d_2)$. Thus the Banach fixed point theorem implies that $T$ has
a unique fixed point. The random variable $X$ has this fixed point as
distribution. The same line of thought actually implies that
$d_2(X_n,X)\to0$. A formal proof can be found in~\cite{NeRu1999}. As
stated in (\ref{Mallow}), the convergence in $(\mathscr M_2^0, d_2)$ is
strong enough to imply convergence of second moments. In particular,
\[
\operatorname{\mathbf{Var}}(\Psi(T^n) )\sim\zeta n^2,
\]
where $\zeta=\operatorname{\mathbf{Var}}(X)$. Computing $\mathbf
{E}[X^2]$ using the fixed point
equation, one easily obtains the following expression for $\zeta$:
%
\begin{equation}\label{equiv}
\zeta=\operatorname{\mathbf{Var}}(X)=\frac{\mu^{-2}\mathbf
{E}[(\sum_{i=1}^{b}V_i\log V_i)^2]-1}{1-\sum
_{i=1}^{b}\mathbf{E}[V_i^2]}.
\end{equation}
This expression may also be obtained using estimates based on renewal
theory in the spirit of our proof of Theorem~\ref{expectedpath}.

\section{Precise asymptotics for the average path length}\label{secaverage}
\subsection{\texorpdfstring{Plan of the proof of Theorem \protect\ref{expectedpath}}
{Plan of the proof of Theorem 3.1}}\label{secplan}

In the previous section, we have explained why precise asymptotics for
$\mathbf{E}[\Psi(T^n)]$ imply convergence in distribution of $\Psi(T^n)$
(suitably rescaled). We now move on to the proof of Theorem~\ref{expectedpath}.

Recall that $D_i$ denotes the depth of the $i$th inserted item. Write
$i\in T_u$ if the item $i$ is stored in the subtree rooted at $u$. Then
rearranging the sum in the definition of~$\Psi(T^n)$, we see that
%
\begin{equation}\label{eqsumsizes}
\Psi(T^n)=\sum_{i=1}^n D_i = \sum_{i=1}^n \sum_{u\ne\sigma} \one
_{ \{ i\in T_u \} }=\sum_{u\ne\sigma} n_u.
\end{equation}

Recall the following fact, which we used already in Section
\ref{seccontraction}:
\[
\frac1 n \operatorname{Mult}(n; V_1,\ldots, V_b) \to(V_1,\ldots, V_b),
\]
almost surely, as $n\to\infty$. We actually have a similar behavior for
any random variable $n_v$, when $v$ is a fixed node (so in particular,
its depth does not depend on~$n$). For a node $u$, the components
$V_1,V_2,\ldots,V_b$ of $\mathcal{V}_u$ are naturally associated to the
children $u_1,u_2,\ldots,u_b$ of $u$, and we can define $V_{u_i}=V_i$.
For the root node~$\varnothing$, define $V_\varnothing=1$. Then let
%
\begin{equation}\label{eqdefLu}
L_u=\prod_{v\preceq u} V_v,
\end{equation}
where $v\preceq u$ if $v$ is an ancestor of $u$. The random variables
$(L_u, u\in\mathcal U)$ define a recursive partition of $[0,1]$, where
$L_u$ is the \textit{length} of the interval associated with~$u$. In
general, for any fixed node $u$, we have
\[
\frac{n_v}n\to L_v,
\]
almost surely as $n\to\infty$. So, as long as $n_v$ is large it should
be well approximated by $n L_v$. This suggests that the sum in (\ref
{eqsumsizes}) be decomposed into the contributions of the top and of
the fringe of the tree. We define the separation in terms of a
parameter $B$ measuring the size of the trees pending in the fringe.
The lengths $L_v$ are decreasing on any path from the root. So let~$R$
be the collection of nodes such that $r\in R$ if $r$ has $nL_r<B$ but
for all its strict ancestors $v$ we have $nL_v\geq B$. We write
$T_{r}, r\in R$, for the subtrees rooted at the nodes that belong to~$R$.

Then
%
\begin{equation}\label{eqdecompositionsum}
\mathbf{E}[\Psi(T^n)]=\mathbf{E}\biggl[\sum_{v\ne\varnothing} n_v\one
_{ \{ nL_v\geq B \} }\biggr]+\mathbf{E}\biggl[\sum_{r\in R}\Psi(T^{n_r})+n_r\biggr],
\end{equation}
since given $n_r$, the total path length of $T_r$, $r\in R$, is
distributed like $T^{n_r}$. [The term $n_r$ needs to be added since the
cardinality of the root of a tree~$T$ is not taken\vadjust{\goodbreak} into account from
our definition of $\Psi(T)$.] The following two propositions gather the
asymptotics for the two terms in (\ref{eqdecompositionsum}) above
that will enable us to prove Theorem~\ref{expectedpath}. In the
following, we let
\[
d=\sup\{a\ge0\dvtx {\mathbf P}(\ln V\in a \mathbb Z)=1\}.
\]
Indeed, as we already mentioned (it will become clear soon), the
arithmetic properties of $\ln V$ influence the asymptotics.
%
\begin{prop}\label{lemma1}
There exists a constant $K$ such that, for all $n$ large enough, and
all $B$, we have
\[
\biggl|\mathbf{E}\biggl[\sum_{v\ne\varnothing} n_v\one_{ \{ n L_v\geq B
\} }\biggr]-\frac1 \mu n
\ln
\biggl(\frac n B\biggr)- n \phi_1\biggl(\ln\frac n B\biggr)\biggr|\le K
\frac n B ,
\]
where $\mu$ is the constant in (\ref{mu}) and $\phi_1$ is a continuous
$d$-periodic function; in particular, $\phi_1$ is constant when $d=0$.
\end{prop}
%
\begin{prop}\label{lemma2}
There exists a constant $K$ such that, for all $n$ large enough, all
$\varepsilon>0$ small enough and $B=\varepsilon^{-8}$, we have
%
\begin{equation}\label{nonlattice}
\biggl|\mathbf{E}\biggl[\sum_{r\in R}\Psi(T^{n_r})+n_r\biggr]- n
\varphi_B \biggl(\ln \frac n B\biggr) \biggr|\le K \varepsilon n
\end{equation}
for some $\varphi_B$, a $d$-periodic function that depends on $B$.
Furthermore, there exists a constant $K'$ (independent of $B$) such
that, for $\varepsilon>0$ small enough,
%
\begin{equation}\label{eqcontinuityphiB}
{\sup_{|q-q'|\le\varepsilon^3}} |\varphi_B(q)-\varphi_B(q')|\le K'
\varepsilon
\ln(1/\varepsilon).
\end{equation}
%
\end{prop}



The proofs of Propositions~\ref{lemma1} and~\ref{lemma2} both rely on
renewal theory: first, the sum $S_{n,B}$ is easily approximated by a
function of sums of i.i.d. random variables; second, the sizes $n_r$
in the second contribution can be estimated using overshoot arguments.
The necessary technical lemmas are introduced in the following section.
Then, we prove Propositions~\ref{lemma1} and~\ref{lemma2} in
Sections~\ref{sectop} and~\ref{secfringe}, respectively.

Before we proceed to the proofs of Propositions~\ref{lemma1} and \ref
{lemma2}, we prove that they indeed imply Theorem~\ref{expectedpath}.
The nonlattice case should be rather clear, but the lattice case
requires a little care.
\begin{pf*}{Proof of Theorem~\ref{expectedpath}}
We have been precise in the statements of Propositions~\ref{lemma1}
and~\ref{lemma2}; we now take the liberty to use $O( \cdot)$
notation to simplify the discussion. It is understood that the hidden
constants do not depend on $n, \varepsilon$ or $B$.

\begin{longlist}
\item
First assume that
$\ln V$ is nonlattice ($d=0$). Let $n,\widehat n$ be integers such
that $n\le\widehat n$. Fix $\varepsilon>0$, and choose $B=\varepsilon^{-20}$.
Then by the triangle inequality and Propositions~\ref{lemma1} and
\ref{lemma2},
%
\[
\biggl|\biggl(\frac{\mathbf{E}[\Psi(T^n)]}{n}-\mu^{-1}\ln n\biggr)-
\biggl(\frac{\mathbf{E}[\Psi(T^{\widehat{n}})]}{\widehat{n}}-\mu
^{-1}\ln
\widehat n\biggr)\biggr|=O(\varepsilon)\vadjust{\goodbreak}
\]
as $n\to\infty$. Thus, the sequence $(n^{-1} \mathbf{E}[\Psi
(T^n)]-\mu
^{-1}\ln
n, n\ge0)$ is Cauchy, hence the result.

\item If $\ln V$ is lattice, the situation is different since we cannot
directly invoke similar arguments. In particular, we need to prove the
existence and continuity of the function $\varpi$. Fix $\beta\in[0,d)$
and consider $\Omega_\beta=\{n\ge1\dvtx\break \exists k\in\mathbb N, |{\ln
n-kd+\beta}|\le n^{-1}\}$, the set of integers such that $\ln n \mbox
{ mod } d$ is close to $\beta$. Then, by the triangle inequality and
Propositions~\ref{lemma1} and~\ref{lemma2}, we have
\begin{eqnarray*}
&&\biggl|\biggl(\frac{\mathbf{E}[\Psi(T^n)]}{n}-\mu^{-1}\ln n\biggr)-
\biggl(\frac{\mathbf{E}[\Psi(T^{\widehat{n}})]}{\widehat{n}}-\mu
^{-1}\ln
\widehat n\biggr)\biggr|\\
&&\qquad\le\biggl|\phi_1\biggl(\ln\frac n B\biggr) - \phi_1\biggl(\ln\frac{\widehat
n}B\biggr)\biggr|+\biggl|\varphi_B\biggl(\ln\frac n B\biggr) - \varphi_B\biggl(\ln\frac
{\widehat n} B\biggr)\biggr|\\
&&\qquad\quad{}+O(\varepsilon)+O(1/B)\\
&&\qquad= |\phi_1(\ln n) - \phi_1(\ln\widehat n)|+
|\varphi
_B(\ln n) - \varphi_B(\ln\widehat n)|+O(\varepsilon),
\end{eqnarray*}
if we choose $\varepsilon$ in such a way that $B=\varepsilon^{-20}=\beta
\mbox
{ mod }d$. Now, $\phi_1$ is continuous and $d$-periodic so that there
exists $n_0$ (independent of $\beta$) such that $|\phi_1(\ln n)-\phi
_1(\ln\widehat n)|\le\varepsilon$ when $n,\widehat n\ge n_0$ inside
$\Omega_\beta$. On the other hand, for $n,\widehat n\in\Omega_\beta$
such that $n,\widehat n\ge2 \varepsilon^{-3}$, we have
\[
|\varphi_B(\ln n)-\varphi_B(\ln\widehat n)|\le K' \varepsilon\ln
(1/\varepsilon).
\]
Note that the bounds obtained are all uniform in $\beta$. It follows
that for every $\varepsilon>0$, there exists $n_1=\max\{n_0,\varepsilon
^{-3}\}
$ such that for $n,\widehat n\in\Omega_\beta$ satisfying $n,\widehat
n\ge n_1$, we have
\[
\biggl|\biggl(\frac{\mathbf{E}[\Psi(T^n)]}{n}-\mu^{-1}\ln n\biggr)-
\biggl(\frac{\mathbf{E}[\Psi(T^{\widehat{n}})]}{\widehat{n}}-\mu
^{-1}\ln
\widehat n\biggr)\biggr|\le O(\varepsilon)+K' \varepsilon\ln(1/\varepsilon).
\]
Therefore, the subsequences $(n^{-1} \mathbf{E}[\Psi(T^n)]-\mu
^{-1}\ln n,
n\in
\Omega_\beta)$, $\beta\in[0,d)$, are uniformly Cauchy (in $\beta
$). It
follows that there exists a fixed function $\varpi$ defined on $[0,d)$
such that, for every $\beta$ and $n \in\Omega_\beta$,
\[
\mathbf{E}[\Psi(T^n)]=\frac1 \mu n \ln n + n \varpi(\beta)+o(n).
\]
Furthermore, the function $\varpi$ is continuous. This is easily seen
using the same arguments with $n\in\Omega_\beta$, $\widehat n\in
\Omega
_{\beta'}$ and $|\beta-\beta|<\varepsilon$. Once the definition of
$\varpi
$ is extended by periodicity, the continuity ensures that we can write
the asymptotics for $\mathbf{E}[\Psi(T^n)]$ in the form claimed in
(\ref
{exptotal1h}). This completes the proof in the lattice case.\qed
\end{longlist}
\noqed\end{pf*}

\subsection{The renewal structure of split trees}
Renewal theory has already been used for studying random trees in
\cite{Holmgren2009,Holmgren2010a,MoRo2005,MoRo2010,Janson2010a}. The
present paper is another example of its wide applicability. We start by
quantifying the deviation between $n_v$ and $nL_v$ for fixed nodes
$v\in\mathcal U$.\vadjust{\goodbreak}
%
\begin{Lemma}\label{cardinality}For any node $v$, we have for all $x$
large enough
\[
\mathbf{P}\bigl(|n_v-n L_v| > (nL_v)^{2/3} \mid nL_v>x \bigr) \le x^{-1/4}.
\]
\end{Lemma}
\begin{pf}
First note that by the triangle inequality
\begin{eqnarray*}
&&\mathbf{P}\bigl(|n_v-nL_v| > (nL_v)^{2/3} \mid nL_v>s \bigr)\\
&&\qquad\le\mathbf{P}\bigl(2|n_v - \bin(n,L_v)|> (nL_v)^{2/3} \mid nL_v>x \bigr)\\
&&\qquad\quad{} +\mathbf{P}\bigl(2|{\bin}(n,L_v)-n L_v|>(nL_v)^{2/3} \mid nL_v>x \bigr).
\end{eqnarray*}
Suppose that $|v|=d$ and let $\mathscr{G}_d$ be the $\sigma$-field
generated by the random variables $V_u$ for $|u|\le d$. Conditioning on
$\mathscr G_{d}$, the recursive splits of the cardinalities $n_v$
defined in (\ref{mult}) give in a stochastic sense the following bound
for~$n_v$:
%
\begin{equation}\label{Binomial2}
|n_v-\operatorname{Bin}(n,L_v)|\leq_{st} \sum_{u\preceq
v}\operatorname{Bin}(s,L_v/L_u).
\end{equation}
Now, by (\ref{Binomial2}), Chebyshev's inequality and Chernoff's bound
for binomials (see, e.g.,
\cite{Chernoff1952,Hoeffding1963,JaLuRu2000}) we obtain
\begin{eqnarray*}
&&\mathbf{P}\bigl(|n_v-nL_v| > (nL_v)^{2/3} \mid nL_v>x \bigr)\\
&&\qquad\leq
2x^{-2/3}\mathbf{E}\biggl[\sum_{u\preceq v}\operatorname{Bin}(s,L_v/L_u)\biggr]\\
&&\qquad\quad{} +
\mathbf{E}\biggl[\exp\biggl(\frac{- (nL_v)^{4/3}}{8(nL_v+(nL_v)^{2/3}/6)}\biggr)
\Bigm|
nL_v>x\biggr]\\
&&\qquad\leq2 s x^{-2/3}\sum_{k\ge0}b^{-k}+e^{-x^{1/4}}
\leq x^{-1/4}
\end{eqnarray*}
for all $x$ large enough.
\end{pf}
%

When the cardinalities $n_v$ are close to the product $n L_v$, renewal
theory allows us to get approximations suitable to
prove Propositions~\ref{lemma1} and~\ref{lemma2}. It is
convenient to introduce the additive form $S_v=-\ln L_v$. For $|v|=k$,
\[
S_v\eqdist S_k=\sum_{i=1}^{|v|} -\ln V_i,
\]
where $V_i$, $i\ge1$, are i.i.d. copies of $V$.
We define the exponential renewal function
%
\begin{equation}\label{renewalfunction}
U(t):=\sum_{k=1}^{\infty}b^k\mathbf{P}(S_k\leq t),
\end{equation}
which satisfies the following renewal equation with $\nu(t)=b\mathbf
{P}(-\ln V\leq t)$:
%
\begin{equation}\label{eqnsrenewal}\qquad
U(t)=\nu(t)+(U\ast d\nu)(t) \qquad\mbox{where } (U\ast d\nu
)(t)=\int_0^t U(t-z) \,d\nu(z).
\end{equation}
The measure $d\nu(t)$ is not a probability measure. To work with more
convenient renewal equations, involving probability measures, we
introduce the tilted measure $d\omega(t)=e^{-t}\,d\nu(t)$. It is easily
seen that $d\omega(t)$ is a probability measure, and defines a random
variable $X$ by ${\mathbf P}(X\in dt)=d\omega(t)$. In fact~$\omega$
is the
distribution function of $-\ln\Delta$, where $\Delta$ is the
size-biased random variable in (\ref{mu}): writing $I$ for a random
variable that is $i$ with probability $V_i$ given $(V_1,\ldots, V_b)$,
we have
\begin{eqnarray*}
{\mathbf P}(-\ln\Delta\leq x)
&=&{\mathbf E}\mathbf{E}\bigl[ \one_{ \{ -\ln V_I\leq x \} } \mid
(V_1,\ldots,V_b)\bigr]\\
&=&\mathbf{E}\Biggl[\sum_{i=1}^{b}\one_{ \{ -\ln V_i\leq x \} }V_i\Biggr]\\
&=&b\mathbf{E}\bigl[\one_{ \{ -\ln V\leq x \} }e^{-\ln V}\bigr]=\omega(x).
\end{eqnarray*}

Then, from (\ref{mu}), $X$ obviously satisfies
\[
\mathbf{E}[X]=\mathbf{E}[-\ln\Delta]=\mu\quad\mbox{and}\quad
\mathbf{E}[X^2]=\sigma^2+\mu^2.
\]
The renewal equation (\ref{eqnsrenewal}) can then be rewritten as
%
\begin{equation}\label{ny} \widehat{U}(t)=\widehat{\nu
}(t)+(\widehat
{U}*d\omega)(t),
\end{equation}
where $\widehat{U}(t):=e^{-t}U(t)$ and $\widehat{\nu}(t):=e^{-t}\nu
(t)$. The first-order asymptotics for $U(t)$ as $t\to\infty$ follows
from the standard renewal theorem applied to $\widehat U(t)$ (see also
Theorem 7.1, Chapter V of~\cite{Asmussen2003} or Lemma 3.1 of
\cite{Holmgren2009} for a formal proof):
%
\begin{equation}\label{renewalequation3}
U(t)=\widehat U(t)e^t=\mu^{-1} e^t +o(e^t),\qquad t\rightarrow\infty.
\end{equation}

We will need some information about the second-order behavior of
$U(t)$. The following lemma will be sufficient for us.
%
\begin{Lemma}\label{Vlem}
Let $d=\sup\{a\ge0\dvtx {\mathbf P}(\ln V\in a \mathbb Z)=1\}$, so that
$d=0$ if
$\ln V$ is nonlattice. Then, as $x\rightarrow\infty$
%
\begin{eqnarray} \label{V}
&&
\int_{0}^{x}e^{-t}\bigl(U(t)-\mu^{-1}e^t\bigr)\,dt\nonumber\\[-8pt]\\[-8pt]
&&\qquad=
\cases{
\displaystyle \frac{\sigma^2-\mu^2}{2\mu^2}-\mu^{-1}+o(1),&\quad if $d=0$,\vspace*{2pt}\cr
\displaystyle \frac{\sigma^2-\mu^2}{2\mu^2}-\mu^{-1}+\phi(x)+o(1), &\quad if
$d>0$,}\nonumber
\end{eqnarray}
%
where $\phi(x)$ is a bounded continuous periodic function with period
$d$. 
\end{Lemma}
\begin{pf}Let $X_k$ be i.i.d. copies of a random variable $X$
defined by ${\mathbf P}(X\in dt)=e^{-t}\,d\nu(t)$.
Define the (standard) renewal function
%
\begin{equation}\label{eqstandardrenewal}
F(t):=\sum_{n\ge0} {\mathbf P}\Biggl(\sum_{k=1}^n X_k\leq t\Biggr).
\end{equation}
Then the renewal theorem (Theorem V.2.4 of~\cite{Asmussen2003}) applied
to (\ref{ny}) yields
%
\begin{equation}\label{standardrenewal4}
e^{-t}U(t)=\widehat{U}(t)=\int_{0}^{t}\widehat{\nu}(t-u)\,dF(u)=\int
_{0}^{\infty}\widehat{\nu}(u)\,dF(t-u).
\end{equation}
[Note that $dF(t)$ includes a term $d{\mathbf P}(0\le t)=\delta_0(t)$.]
By Fubini's theorem we obtain
%
\begin{eqnarray}\label{standardrenewal5}\int_{0}^{x}e^{-t}\bigl(U(t)-\mu
^{-1}e^t\bigr)\,dt&=&
\int_{0}^{\infty}\widehat{\nu}(u)\int_{0}^{x}dF(t-u)\,du-\frac{
x}{\mu
}\nonumber\\[-8pt]\\[-8pt]
&=&\int_{0}^{\infty}\widehat{\nu}(u)F(x-u)\,du-\frac{x}{\mu
}.\nonumber
\end{eqnarray}
Recall that $\widehat{\nu}(x)=\nu(x)e^{-x}$. Integration by parts gives
%
\begin{equation}\label{reneq}
\int_{0}^{\infty}\widehat{\nu}(x)\,dx=b[-e^{-t}{\mathbf P}(-\ln
V\leq t)
]_{0}^{\infty}+\int_{0}^{\infty}e^{-t}\,d\nu(t)=b\mathbf{E}[e^{-\ln
V}]=1.\hspace*{-28pt}
\end{equation}
Rewriting (\ref{standardrenewal5}) as a single integral, it follows that
%
\begin{eqnarray}\label{standardrenewal6}
&& \int_{0}^{x}e^{-t}\bigl(U(t)-\mu
^{-1}e^t\bigr)\,dt\nonumber\\
&&\qquad=\int_{0}^{\infty}\widehat{\nu}(u)\biggl(F(x-u)-\frac{x}{\mu
}
\biggr)\,du\nonumber\\[-8pt]\\[-8pt]
&&\qquad=-\frac{1}{\mu}\int_{0}^{x}\widehat{\nu}(u)u\,du-
\frac{1}{\mu}\int_{x}^{\infty}\widehat{\nu}(u)x\,du\nonumber\\
&&\qquad\quad{}+\int
_{0}^{x}\widehat
{\nu}(u)\biggl(F(x-u)-\frac{x-u}{\mu}\biggr)\,du.\nonumber
\end{eqnarray}
%
We start with the first two terms in (\ref{standardrenewal6}). Using
again integration by parts and applying (\ref{reneq}) yields
%
\begin{eqnarray}\label{parts}
\int_{0}^{\infty}\widehat{\nu}(u)u\,du
&=&\int_{0}^{\infty}e^{-u}{\nu}(u)u\,du\nonumber\\
&=&\int_{0}^{\infty}\widehat{\nu}(u)\,du+\int_{0}^{\infty}ue^{-u}\,
d{\nu
}(u)\\
&=&1+b\mathbf{E}[-V\ln V]=1+\mu,\nonumber
\end{eqnarray}
where the last equality follows from the definition of $\mu$ in (\ref
{mu}). Finally, note that for all $x$,
%
\begin{equation}\label{eqsecterm}
\int_{x}^{\infty}\widehat{\nu}(u)x\,du\le\int_x^\infty\widehat\nu
(u) u\,
du \to0
\end{equation}
as $x\to\infty$ since $\int_{0}^{\infty}|\widehat{\nu
}(u)u|\,du<\infty$.

So it only remains to estimate the third term in (\ref
{standardrenewal6}). This is related to the asymptotics for the renewal
function $F(t)$, which are different depending on whether $\ln V$ is
lattice or not. Write $\{x\}$ for the fractional part of a real number
$x$, that is, $\{x\}=x-\lfloor x\rfloor$. Then, by Theorem 5.1 in
\cite{Gut2009} we have, as $t\to\infty$,
\[
F(t)-\frac{t}{\mu}=\frac{\sigma^2+\mu^2}{2\mu^2}+o(1)
\quad\mbox{and}\quad
F(t)-\frac{t}{\mu}= \frac{\sigma^2+\mu^2}{2\mu^2}+\frac
{d}{\mu
}\biggl(\frac{1}{2}-\biggl\{\frac{t}{d}\biggr\}\biggr)+o(1)
\]
in the nonlattice and the $d$-lattice case, respectively.
Furthermore, by Lorden's inequality (\cite{Lorden1970}, Theorem 1),
\[
0\leq F(t)-\frac{t}{\mu}\leq\frac{\sigma^2+\mu^2}{\mu^2}.
\]

(i) We now first assume that $\ln V$ is nonlattice. The dominated
convergence theorem applied to the last integral in (\ref
{standardrenewal6}), and (\ref{reneq}), yield
%
\begin{eqnarray}\label{standardrenewal7}
\lim_{x\rightarrow\infty} \int_{0}^{\infty}\widehat{\nu}(u)
\biggl(F(x-u)-\frac{x-u}{\mu}\biggr)\one_{ \{ u\leq x \} }\,du&=&\int
_{0}^{\infty
}\widehat{\nu
}(u)\frac{\sigma^2+\mu^2}{2\mu^2}\,du\nonumber\\[-8pt]\\[-8pt]
&=&\frac{\sigma^2+\mu^2}{2\mu^2}.\nonumber
\end{eqnarray}
Putting (\ref{standardrenewal7}) together with (\ref
{standardrenewal6}), (\ref{parts}) and (\ref{eqsecterm}) we obtain, as
$x\to\infty$,
\[
\int_{0}^{x}e^{-t}\bigl(U(t)-\mu^{-1}e^t\bigr)\,dt=-\frac1 \mu-1 + \frac
{\sigma^2
+ \mu^2}{2\mu^2}+o(1),
\]
which proves the claim in (\ref{V}) in the nonlattice case.

(ii) Similarly in the lattice case with span $d$, from the dominated
convergence theorem we obtain
%
\begin{eqnarray}\label{standardrenewal7,1}
&&\int_{0}^{x}\widehat{\nu}(u)\biggl(F(x-u)-\frac{x-u}{\mu}\biggr)\,du\nonumber\\
&&\qquad=\frac{\sigma^2+\mu^2}{2\mu^2}+\frac{d}{\mu} \int_0^x \biggl(\frac
{1}{2}-\biggl\{\frac{x-u}{d}\biggr\}\biggr) \widehat\nu(u)\,du +o(1)\\
&&\qquad=\frac{\sigma^2+\mu^2}{2\mu^2}+\frac{d}{\mu} \int_0^\infty
\biggl(\frac{1}{2}-\biggl\{\frac{x-u}{d}\biggr\}\biggr) \widehat\nu
(u)\,du+o(1)\nonumber
\end{eqnarray}
by (\ref{eqsecterm}).
The function $\phi$ defined for $x\ge0$ by
\[
\phi(x)=\frac{d}{\mu} \int_0^\infty\biggl(\frac{1}{2}-\biggl\{ \frac
{x-u}{d}\biggr\}\biggr) \widehat\nu(u)\,du\vadjust{\goodbreak}
\]
is clearly $d$-periodic. Furthermore, the function $\phi( \cdot)$ is
continuous. Indeed, for any $x,y$ such that $|x-y|< \varepsilon$ we have
\begin{eqnarray*}
\phi(y)
&=&\frac{d}{\mu} \int_0^\infty\biggl(\frac{1}{2}-\biggl\{\frac
{y-u}{d}\biggr\}\biggr) \widehat\nu(u)\,du\\
&=&\frac{d}{\mu} \int_0^\infty\biggl(\frac{1}{2}-\biggl\{\frac
{y-u}{d}\biggr\}\biggr) \one_{ \{ y-u\ \mathrm{mod}\ d \in[\varepsilon
,1-\varepsilon] \} }
\widehat\nu(u)\,du\\
&&{}
+ \frac{d}{\mu} \int_0^\infty\biggl(\frac{1}{2}-\biggl\{\frac
{y-u}{d}\biggr\}\biggr) \one_{ \{ y-u\ \mathrm{mod}\ d \notin[\varepsilon
,1-\varepsilon] \} }
\widehat\nu(u)\,du.
\end{eqnarray*}
It follows that
\begin{eqnarray*}
|\phi(y)-\phi(x)|
&\le&\frac2 \mu\varepsilon+ 2 \sup_{z\in\{x,y\}} \frac{d}{\mu}
\int
_0^\infty\biggl|\frac{1}{2}-\biggl\{\frac{z-u}{d}\biggr\}\biggr|
\one_{ \{ z-u \ \mathrm{mod}\ d \notin[\varepsilon,1-\varepsilon] \} }
\widehat\nu(u)\,du\\
&\le&\frac2 \mu\varepsilon+2 \sup_{z\in\{x,y\}} \frac d \mu\int
_0^\infty\one_{ \{ z-u\ \mathrm{mod}\ d \notin[\varepsilon,1-\varepsilon
] \} }
\widehat
\nu(u)\,du.
\end{eqnarray*}
Since $|\widehat\nu(u)|=e^{-u}b {\mathbf P}(-\ln V\le t)\le b$, the dominated
convergence theorem implies that $|\phi(y)-\phi(x)|\to0$ as
$\varepsilon
\to0$.

Finally, putting (\ref{standardrenewal7,1}) together with (\ref
{standardrenewal6}), (\ref{parts}) and (\ref{eqsecterm}) as before
proves the lattice case in (\ref{V}).
\end{pf}

\subsection{Contribution of the top of the tree}\label{sectop}

In this section, we prove Proposition~\ref{lemma1}. For the top of the
tree, the sizes $n_v$ are well approximated by $\bin(n, L_v)$. This
suggests that the main contribution of the top of the tree should be
%
\begin{equation}\label{eqdefRnb}
\mathbf{E}\biggl[\sum_{v\ne\varnothing} n_v \one_{ \{ nL_v\ge B \}
}\biggr]=\mathbf{E}\biggl[\sum_{v\ne\varnothing} \bin(n,L_v)\one_{ \{
nL_v\ge B \} }\biggr]+R_{n,B}
\end{equation}
for a remainder $R_{n,B}$ that should be small. We first estimate the
main contribution; we will then quantify $R_{n,B}$ using (\ref{Binomial2}).
%
\begin{Lemma}\label{lemma3}
Let $d=\sup\{a\dvtx{\mathbf P}(\ln V \in a \mathbb Z)=1\}$, so that $d=0$
if $\ln
V$ is nonlattice.
Then, as $n/B\to\infty$,
\begin{eqnarray*}
&&\mathbf{E}\biggl[\sum_{v\ne\varnothing} \bin(n,L_v)\one_{ \{ nL_v\geq B
\} }\biggr]
\\
&&\qquad=
\cases{
\displaystyle \frac1 \mu n \ln\biggl(\frac nB\biggr)+ n
\frac{\sigma^2-\mu^2}{2\mu^2}+o(n), &\quad if $d=0$,\vspace*{2pt}\cr
\displaystyle \frac1 \mu n \ln\biggl(\frac nB\biggr)+ n
\frac{\sigma^2-\mu^2}{2\mu ^2}+n\phi \biggl(\ln\frac n B\biggr)+
o(n),&\quad if $d>0$,}
\end{eqnarray*}
%
where $\mu$ and $\sigma$ are the constants in (\ref{mu}) and $\phi
(\cdot
)$ is a bounded continuous $d$-periodic function.
\end{Lemma}
\begin{pf}
Let $V_i$, $i\ge1$ be i.i.d. copies of $V$, and define $L_k=\prod
_{i=1}^k V_i$ and $S_k=-\ln L_k$. Then, we have
\begin{eqnarray*}
\mathbf{E}\biggl[\sum_{v\ne\varnothing} \bin(n,L_v)\one_{ \{ nL_v\geq B
\} }\biggr]
&=& n \mathbf{E}\biggl[\sum_{k\ge1} b^k L_k \one_{ \{ n L_k\geq B \} }\biggr]\\
&=& n \mathbf{E}\biggl[\sum_{k\ge1} b^d e^{-S_k}\one_{ \{ S_k\leq\ln n
-\ln B \} }\biggr]\\
&=& n \int_0^{\ln(n/B)}\sum_{k\ge1} b^k e^{-t}\,d{\mathbf P}(S_k\leq
t)\\
&=& n \int_0^{\ln(n/B)} e^{-t} \,dU(t),
\end{eqnarray*}
where $U(t)$ is the renewal function defined in (\ref{renewalfunction}).
Using integration by parts we obtain, if $-\ln V$ is nonlattice,
\begin{eqnarray*}
&&
\int^{\ln(n/B)}_0e^{-t}\,dU(t)
\\
&&\qquad=[e^{-t}U(t)]^{\ln(n/B)}_0+\int^{\ln
(n/B)}_0e^{-t}U(t)\,dt\\
&&\qquad=\frac B n U\bigl(\ln(n/B)\bigr)+\int^{\ln(n/B)}_0e^{-t}\bigl(U(t)-\mu
^{-1}e^t\bigr)\,dt\\
&&\qquad\quad{}+\mu
^{-1} \ln(n/B)\\
&&\qquad=\mu^{-1} +o(1) + \frac{\sigma^2-\mu^2}{2\mu^2}-\mu^{-1}\\
&&\qquad\quad{}+\mu
^{-1} \ln
(n/B) +o(1)
\end{eqnarray*}
by Lemma~\ref{Vlem} and (\ref{renewalequation3}). Similarly if $-\ln
V$ is lattice with span $d$, Lemma~\ref{Vlem} and (\ref{renewalequation3}) yield
\begin{eqnarray*}
\int^{\ln(n/B)}_0e^{-t}\,dU(t)
&=&\mu^{-1} +o(1) + \frac{\sigma^2-\mu^2}{2\mu^2}-\mu^{-1}+\mu
^{-1} \ln
(n/B) \\
&&{}+\phi\bigl(\ln(n/B)\bigr) +o(1),
\end{eqnarray*}
where $\phi(t)$ is a continuous periodic function with period $d$. 
%
\end{pf}

We now deal with the remainder $R_{n,B}$ introduced in (\ref
{eqdefRnb}). The difference between $n_v$ and the binomial is bounded
in (\ref{Binomial2}) and we have
\[
|R_{n,B}|\le\mathbf{E}\biggl[\sum_{v\ne\varnothing} \one_{ \{ n L_v\ge
B \} } \sum_{u\preceq v}\bin(s, L_v/L_u)\biggr].
\]

\begin{Lemma}\label{lemma4}
The following estimate holds: there exist a constant and $n_0$ such
that, for every fixed $B$ and $n\ge n_0$, we have
\[
\mathbf{E}\biggl[\sum_{v\ne\varnothing} \one_{ \{ n L_v\ge B \} } \sum
_{u\preceq v}\bin(s, L_v/L_u)\biggr]=O \biggl(\frac{n}{B}\biggr).
\]
\end{Lemma}
\begin{pf}
In the following, $|v|=d$, $|u|=k\le d$, and we write $\ell=d-k$. Then
$L_v$ is distributed as $L_d=L_k \cdot L_\ell$, where the two factors
are products of $k$ and $\ell$ copies of $V$, respectively; all of them
are independent. Swapping the sums over $u$ and $v$, we obtain
%
\begin{eqnarray}\label{salt}
&&\mathbf{E}\biggl[\sum_{v\ne\varnothing} \one_{ \{ nL_v\geq B \} } \sum
_{u\preceq v}\bin(s,L_v/L_u)\biggr]
\nonumber\\
&&\qquad=\mathbf{E}\biggl[\sum_u \sum_{v\dvtx u\preceq v, v\ne\varnothing} s \frac
{L_v}{L_u} \one_{ \{ nL_v\ge B \} }\biggr]
\le s \mathbf{E}\biggl[\sum_{k\ge0} b^k \sum_{\ell\ge0} b^\ell
L_\ell\one_{ \{ n L_k L_\ell\ge B \} }\biggr]\\
%
&&\qquad=s \mathbf{E}\biggl[\sum_{k\ge0}b^k\sum_{\ell\ge0}b^\ell e^{-S_\ell
}\one_{ \{ e^{S_k+S_\ell}\leq n/B \} }\biggr].\nonumber
\end{eqnarray}
First conditioning on $S_k$ in each term of the sum above, and
recalling the renewal function $U(t)$ defined in (\ref
{renewalfunction}), we see that
\[
\mathbf{E}\biggl[\sum_{\ell\ge0}b^\ell e^{-S_\ell}\one_{ \{
e^{S_k+S_\ell}\leq n/B \} } \Bigm| S_k\biggr]
=\int^{\ln(n/B)-S_k}_0 e^{-t}\,dU(t)+b\one_{ \{ e^{S_k}\le n/B \} }.
\]
However, there exists a constant $C$ such that, for any real number $x$,
\[
\int_0^x e^{-t} \,dU(t)\le C x\one_{ \{ x\ge0 \} }.
\]
Going back to (\ref{salt}) and choosing $x=\ln(n/B)-S_k$, it follows that
\begin{eqnarray*}
&&
\mathbf{E}\biggl[\sum_{v\ne\varnothing} \one_{ \{ nL_v\geq B \} } \sum
_{u\preceq v}\bin(s,L_v/L_u)\biggr]\\
&&\qquad\leq C \mathbf{E}\biggl[\sum_{k\ge0} b^k \bigl(\ln(n/B)-S_k+b
\bigr)\one_{ \{ S_k\leq\ln(n/B) \} }\biggr]\\
&&\qquad=C \int^{\ln(n/B)}_0 \bigl(\ln(n/B)-t+b\bigr)\,dU(t)\\
&&\qquad=C \bigl[\bigl(\ln(n/B)-t\bigr)U(t)\bigr]^{\ln(n/B)}_0\\
&&\qquad\quad{}+C'
\int^{\ln(n/B)}_0 U(t)\,dt,
\end{eqnarray*}
where the last line follows by integration by parts and we wrote \mbox{$C'=C(1+b)$}.
The claim then follows from (\ref{renewalequation3}).
\end{pf}


\subsection{\texorpdfstring{Contribution of the fringe: Proof of Proposition \protect\ref{lemma2}}
{Contribution of the fringe: Proof of Proposition 4.2}}\label{secfringe}

Finally, we prove Proposition~\ref{lemma2} that deals with the
contribution of the fringe of the tree. Recall that from (\ref
{eqdecompositionsum}), we have to estimate
%
\begin{equation}\label{eqfringereminder}
\mathbf{E}\biggl[\sum_{r\in R}\widetilde\Psi(T^{n_r})\biggr]:=\mathbf{E}\biggl[\sum
_{r\in R} \Psi(T^{n_r})+n_r\biggr],
\end{equation}
where, for convenience, we introduced $\widetilde\Psi(T^k):= \Psi
(T^k)+k$. The proofs here get quite technical at times, and the reader
should bear in mind that we will essentially express the expected value
in (\ref{eqfringereminder}) as a mixture of the expected values of
$\mathbf{E}[\widetilde\Psi(T^k)]$, for $k$ lower than $B$.

For a node $r$, define the conditional expectation $\Gamma_r=\mathbf
{E}[ \widetilde\Psi(T^{n_r}) \mid n_r]$.
First, the first asymptotic order of the expected total path length
implies that
%
\begin{equation}\label{sumi2,2}
\Gamma_{r}=O(n_r\ln n_r).
\end{equation}
%
The next lemma is used to get an error bound for the sum of the
expected total path lengths of the subtrees $T_r, r\in R$,
with
cardinalities $n_r$ that differ from $n L_r$ by at least $B^{2/3}$
items, so that
we only have to bother about the subtrees $T_r, r\in R$, with
cardinalities $n_r$ that are close to $nL_r$.
%
\begin{Lemma}\label{hoaloa}The following error bound holds:
\[
\mathbf{E}\biggl[\sum_{r\in R} n_{r}\ln n_r\one_{ \{ |n_{r}-nL_{r}|\geq
B^{2/3} \} }\biggr]=\mathcal
{O}\biggl(\frac{n\ln B}{B^{1/4}}\biggr).
\]
\end{Lemma}

We omit the proof; it follows by a simple modification of the proof of
Lem\-ma~4.3 of~\cite{Holmgren2009}.
By Lemma~\ref{hoaloa}, we have
\[
\mathbf{E}\biggl[\sum_{r\in R}\widetilde\Psi(T^{n_r})\biggr]=\mathbf{E}\biggl[\sum
_{r\in R} \Gamma_{r}\one_{ \{ |n_{r}-nL_{r}|\leq B^{2/3} \}
}\biggr]+O\biggl(\frac{n\ln B}{B^{1/4}}\biggr).
\]
Define $R'\subseteq R$ to be the set of ``good'' nodes in $R$:
%
\begin{equation} \label{epsilonigen}
R':=\{r\in R\dvtx|n_{r}-n L_{r}|\leq B^{2/3}\}
\end{equation}
and let $R''\subseteq R'$ be the subset of nodes $r\in R'$ that also
satisfy $nL_r> \varepsilon^2$.

We will now explain that it is enough to consider the nodes $r\in R''$.
The approximation of $U(t)$ in (\ref{renewalequation3}) implies that
the expected number of nodes $v$ such that $n L_v\geq B$ is ${O}(n/B)$;
thus, since each node has at most $b$ children,
%
\begin{equation}\label{eqsizeR}
\mathbf{E}[|R|]={O}(n/B)
\end{equation}
as well. Hence, it follows from (\ref{epsilonigen}) that the expected
number of nodes in the $T_{r}$, $r\in R'$, with $n L_r\leq\varepsilon
^{2}B$ is bounded by ${O}(\varepsilon^{2} n)$.
Using this fact yields
%
\begin{equation}\label{limsuphajul}
\mathbf{E}\biggl[\sum_{r\in R}\widetilde\Psi(T^{n_r})\biggr]=\mathbf{E}\biggl[\sum
_{r\in R''} \Gamma_{r}\biggr]+
\mathcal{O}(\varepsilon^{2}n\ln B)+{O}\biggl(\frac{n\ln B}{B^{1/4}}\biggr).
\end{equation}

Because of the concentration of $n_r$ around $nL_r$, the cardinalities
$n_r$ of the nodes $r\in R$ are naturally related to the behavior of
the ``overshoot'' of the renewal process $(-\ln L_k, k\ge0)$, when it
crosses the line $\ln(n/B)$. Estimating the empirical distribution of
the cardinalities of the nodes $r\in R$ will allow us to approximate
the right-hand side above.
So we further subdivide the nodes
$r\in R$ into smaller classes according to the values of $n L_r$, $r\in R$.

Let $ Z=\{B,B-\gamma B,B-2\gamma B,\ldots,\varepsilon^2 B\}$, where we let
$\gamma=\varepsilon^3$.
We write $R_z\subseteq R, z\in Z$,
for the set of nodes $r\in R$, such that $n L_r\in[z-\gamma B,z)$.
Then~(\ref{limsuphajul}) can be rewritten as
%
\begin{equation}\label{limsupnujul}
\mathbf{E}\biggl[\sum_{r\in R}\widetilde\Psi(T^{n_r})\biggr]
=\mathbf{E}\biggl[\sum_{z\in Z}\sum_{r\in R'\cap R_{z}}\Gamma_{r}\biggr]+
{O}(\varepsilon^{2} n\ln B)+{O}\biggl(\frac{n\ln B}{B^{1/4}}\biggr).
\end{equation}
Even in a fixed class $R_z$, not all the nodes have the same
cardinality $n_r$. So, in order to estimate the expected value in (\ref
{limsupnujul}) we need the following lemma that quantifies the
discrepancy of $\mathbf{E}[\Psi(T^n)]$ under small variations of $n$.
%
\begin{Lemma}\label{lemmatot}There exists a constant $C$ such that, for
any natural numbers $n$ and $K$, we have
\[
|\mathbf{E}[\widetilde\Psi{(T^{n+K})}]-\mathbf{E}[\widetilde
\Psi{(T^{n})}]
|\le C K\ln(n+K).
\]
\end{Lemma}
\begin{pf}From the iterative construction, we clearly have $\mathbf
{E}[\widetilde\Psi(T^{n+K})]\ge\mathbf{E}[\widetilde\Psi(T^n)]$;
so it
suffices to bound the increase in path length when adding~$K$ extra
items to the tree $T^n$. Thinking again of the iterative construction,
every ball trickles down until it finds a leaf. Then, either it sits
there if there is room left, or it triggers a growth of the tree. It is
important to notice that only these $s+1$ balls may move. Furthermore,
the increase in depth of any of the $s+1$ items (the last one, plus the
$s$ that were already sitting at the leaf) is at most the height of the
final tree $H_{n+K}$. Hence, upon adding $K$ items, the path length
increases by $K(s+1) H_{n+K}\le C K \ln(n+K)$,
by the results of~\cite{Devroye1998} on the height of split trees.
\end{pf}

Write $f_x=\mathbf{E}[\widetilde\Psi(T^{\lfloor x\rfloor})]$. Then
Lemma \ref
{lemmatot} ensures that, for any node $r\in R'\cap R_z$, we have
$\Gamma_r=f_z + O(\gamma B \ln B)$.
By using (\ref{epsilonigen}) and Lemma~\ref{lemmatot}, from~(\ref
{limsupnujul}) we obtain
%
\begin{eqnarray}\label{limsupnuigen}
\mathbf{E}\biggl[\sum_{r\in R}\widetilde\Psi(T^{n_r})\biggr]
&=&\sum_{z\in Z}\mathbf{E}[|R'\cap R_z|]\bigl(f_z+{O}(\gamma B\ln
B)\bigr)\nonumber\\
&&{}+{O}(\varepsilon
^{2}n \ln B )+{O}\biggl(\frac{n \ln
B}{B^{1/4}}\biggr)\nonumber\\[-8pt]\\[-8pt]
&=&\sum_{z\in Z}\mathbf{E}[|R'\cap R_z|] f_z + O (\gamma n \ln B)
\nonumber\\
&&{}+{O}(\varepsilon
^{2}n \ln B )+{O}\biggl(\frac{n \ln B}{B^{1/4}}\biggr),\nonumber
\end{eqnarray}
since $\mathbf{E}[|R|]= O(n/B)$ by (\ref{eqsizeR}).

So the contribution of the fringe is essentially a mixture of the
$f_z$, $z\in Z$. To complete the proof of Proposition~\ref{lemma2}, it
suffices to estimate the mixing measure $\mathbf{E}[|R'\cap R_z|]$,
$z\in Z$.
We first focus on the asymptotics for $\mathbf{E}[|R_z|]$, $z\in Z$. The
following result is obtained by an application of the key renewal theorem.
%
\begin{Lemma}\label{lem12}
Fix $\varepsilon>0$ and let $S:=\{1, 1-\gamma,1-2\gamma,\ldots
,\varepsilon
^{2}\}$, where \mbox{$\gamma=\varepsilon^3$}.
Let $d=\sup\{a\dvtx {\mathbf P}(\ln V\in a \mathbb Z)=1\}$. If $d>0$, we suppose
that $\ln B\in d\mathbb N$. Then for any $\alpha\in S$ we have, as
$n\to
\infty$,
%
\begin{equation}\label{residual42}\quad
\frac{\mathbf{E}[|R_{\alpha B}|]}{n/B}=
\cases{
c_{\alpha} +o(1), &\quad if $\ln V$ is nonlattice $(d=0)$,\cr
\psi_{\alpha}(\ln n) +o(1), &\quad if $\ln V$ is $d$-lattice $(d>0)$,}
\end{equation}
for a constant $c_{\alpha}$ (only depending on $\alpha$ and $\gamma$),
$\psi_{\alpha}(\cdot)$ is the $d$-periodic function given in (\ref
{hopp2,1}) below.
\end{Lemma}
\begin{pf}Let $V_j, j\ge1$, be i.i.d. copies of $V$. For an
integer $k$, write $S_k=-\sum_{j=1}^{k}\ln V_j$. 
Then, by definition, for $\alpha\in S$, we have
\begin{eqnarray*}
\mathbf{E}[|R_{\alpha B}|]
&=&\sum_{u\in U} {\mathbf P}(u\in R_{\alpha B})\\
&=&\sum_{k=0}^{\infty}b^{k+1} \bigl({\mathbf P}\bigl(S_k-\ln V_{k+1}> \ln
(n/B)-\ln\alpha\mbox{ and }S_k\leq\ln(n/B)\bigr)\\
&&\hspace*{40pt}{} -{\mathbf P}\bigl(S_k- \ln V_{k+1}>\ln(n/B)-\ln
(\alpha-\gamma) \\
&&\hspace*{188pt}\mbox{ and }S_k\leq\ln(n/B)\bigr)\bigr)\\
&=&\int_0^{\ln(n/B)}b{\mathbf P}\bigl(\ln(n/B)-t-\ln\alpha\\
&&\hspace*{57pt}<- \ln
V_{k+1}\leq\ln(n/B)-t-\ln( \alpha-\gamma)\bigr)\,dU_0(t),
\end{eqnarray*}
where $U_0(t)=U(t)+1$ is a simple modification of the renewal
$U(t)=\break\sum_{k\ge1} b^k {\mathbf P}(S_k\le t)$ defined in (\ref
{renewalfunction}). Thus,
seeing $\mathbf{E}[|R_{\alpha B}|]$ as a function of $\ln(n/B)$ and writing
%
\begin{equation}
\label{hopp}\quad
H(q):=\int_{0}^{q}b{\mathbf P}\bigl(q-t-\ln\alpha<- \ln V_{k+1}\leq
q-t-\ln( \alpha-\gamma)\bigr)\,dU_0(t),\vadjust{\goodbreak}
\end{equation}
we have $\mathbf{E}[|R_{\alpha B}|]=H(\ln(n/B))$. So we are after the
asymptotics for $H(q)$, as $q\to\infty$.
%
It is convenient to use a change of measure to relate $H(q)$ to\vadjust{\goodbreak}
a~renewal function associated to a \textit{probability} measure. We have
%
\begin{eqnarray}\label{hopp1}\quad
\widehat H(q)
:\!&=&e^{-q} H(q) \nonumber\\
&=& \int_0^q e^{-(q-t)} G(q-t) e^{-t} \,dU_0(t)\\
&=& \int_0^q b e^{-(q-t)} {\mathbf P}\bigl(q-t-\ln\alpha<- \ln
V_{k+1}\leq q-t-\ln( \alpha-\gamma)\bigr) \,dF(t),\nonumber
\end{eqnarray}
where $F(t)$ is the standard renewal function already introduced in
(\ref{eqstandardrenewal}).
The asymptotics for the integral above are then easily obtained by
using the key renewal theorem. In particular, they depend on whether
$\ln V$ is lattice or not.\vspace*{8pt}

(i) If $\ln V$ is nonlattice,
by the key renewal theorem (\cite{Gut2009}, Theorem II.4.3), we obtain
%
\begin{equation}\label{hopp2}\qquad
\lim_{q\to\infty}\widehat{H}(q)= c_\alpha:=\frac{b}{\mu}\int
_{0}^{\infty}e^{-t}
{\mathbf P}\bigl(t-\ln\alpha<- \ln V\leq t-\ln( \alpha-\gamma)\bigr)\,dt.
\end{equation}
Note that the constant $c_{\alpha}$ only depends on $\alpha$ (and
$\gamma$) and that $\sum_{\alpha\in S}c_\alpha\le b/\mu$.
Thus, since $\widehat{H}(x)=e^{-x}H(x)$ it follows immediately that
$\mathbf{E}[|R_{\alpha B}|]
=\frac{n}{B}c_{\alpha}+o(\frac{n}{B})$ which proves the
nonlattice case in (\ref{residual42}).

(ii) Similarly, if $\ln V$ is lattice with span $d$, the key renewal
theorem (see~\cite{Gut2009}, Theorem II.4.3, or
\cite{Janson2010a}, Theorem A.7) implies that
%
\begin{eqnarray}\label{hopp2,1}\qquad
\widehat{H}(q)&\sim&\psi_{\alpha}(q)\nonumber\\[-8pt]\\[-8pt]
:\!&=& \frac{bd}{\mu}\sum_{k\dvtx
kd\le q} e^{kd-q}
{\mathbf P}\bigl(q-kd-\ln\alpha<- \ln V\leq q-kd-\ln(
\alpha-\gamma)\bigr)\nonumber
\end{eqnarray}
as $q\to\infty$.
Note that $\psi_{\alpha}$ is a (positive) $d$-periodic function.
Observe also that for fixed $\alpha$, the function $\psi_\alpha(
\cdot
)$ is not continuous since $\ln V\in d\mathbb Z$ almost surely.
Since $\widehat{H}(x)=e^{-x}H(x)$, it follows from (\ref{hopp2,1}) that
$\mathbf{E}[|R_{\alpha B}|]
\sim\frac{n}{B}\psi_{\alpha}(\ln(n/B))$. This proves the lattice case
in (\ref{residual42}), and completes the proof.
%
%
\end{pf}

With Lemma~\ref{lem12} in hand, we can now deduce the asymptotics for
$\mathbf{E}[|R'\cap R_z|]$, $z\in Z$ and use them in (\ref
{limsupnuigen}) to
complete the proof of Proposition~\ref{lemma2}. Recall that $R'=\{r\in
R\dvtx |n_{r}-n L_{r}|\leq B^{2/3}\}$.
Clearly, $\mathbf{E}[|R' \cap R_{\alpha B}|]\le\mathbf{E}[|R_{\alpha
B}|]$. Furthermore,
\begin{eqnarray*}
\mathbf{E}[|R'\cap R_{\alpha B}|]
&=&\sum_{r\in R} {\mathbf P}\bigl(|n_r - nL_r|\le B^{2/3}, (\alpha-\gamma
) B\le nL_r< \alpha B\bigr)\\
&=& \sum_{r\in R} {\mathbf P}\bigl((\alpha-\gamma) B \le n L_r < \alpha B\bigr)
\\
&&\hspace*{14pt}{}\times\mathbf{P}\bigl(|n_r - nL_r|\le B^{2/3} \mid(\alpha-\gamma)B\le n L_r <
\alpha B \bigr)\\
&\ge&\mathbf{E}[|R_{\alpha B}|]\bigl(1-O(B^{-1/4})\bigr)
\end{eqnarray*}
by Lemma~\ref{cardinality}.
We now choose $B=\varepsilon^{-20}$ so that
$B^{-1/4}=\varepsilon^{5}$.\vspace*{8pt}

\mbox{}\hphantom{i}(i) If $\ln V$ is nonlattice, it follows from Lemma~\ref{lem12} that
for each choice of~$\gamma$ there is a constant $K_\gamma$ such that
for all $\alpha\in S$ and some constant $c_\alpha$ (that of
Lem\-ma~\ref{lem12}) we have
\[
\biggl|\frac{\mathbf{E}[|R'\cap R_{\alpha B}|]}{n/B}-c_{\alpha}
\biggr|\leq\gamma^{2}+{O}(B^{-1/4})=\gamma^2+O(\varepsilon^5)=O(\varepsilon^5),
\]
whenever $n/B\geq K_{\gamma}$. So for all $n$ large enough, since
$f_x=O(x\ln x)$, we have
\begin{eqnarray*}
\mathbf{E}\biggl[\sum_{r\in R}\widetilde\Psi(T^{n_r})\biggr]
&=&\sum_{\alpha\in S}c_{\alpha}\frac{n}{B}f_{\alpha B}+\frac
{n}{B}\sum
_{\alpha\in S}{O}(f_{\alpha B}\varepsilon^{5})+{O}(n \gamma\ln
B)+{O}(\varepsilon^{2} n \ln B )
\\
&=&n\sum_{\alpha\in S} \frac{f_{\alpha B}}{B}c_{\alpha
}+{O}(\varepsilon n).
\end{eqnarray*}
This proves Proposition~\ref{lemma2} when $\ln V$ is nonlattice.

(ii) Similarly, if $\ln V$ is $d$-lattice, for any choice of $\gamma$,
there is a $K_{\gamma}$ such that for any $\alpha\in S$ and some
continuous $d$-periodic function $\psi_{\alpha}(t)$ [that of
Lemma~\ref{lem12} defined in (\ref{hopp2,1})], we have
\[
\biggl|\frac{\mathbf{E}[|R'\cap R_{\alpha B}|]}{n/B}-\psi_{\alpha
}(\ln
n)\biggr|
\leq\gamma^{2}+{O}(B^{-1/4})=\gamma^2 +O(\varepsilon^5),
\]
whenever $n/B\geq K_{\gamma}$. It follows that
%
%
\begin{eqnarray}\label{eqexprfixedeps}
\mathbf{E}\biggl[\sum_{r\in R}\widetilde\Psi(T^{n_r})\biggr]
&=&\sum_{\alpha\in S}\psi_{\alpha}(\ln n)\frac{n}{B}f_{\alpha
B}+\frac
{n}{B}\sum_{\alpha\in S}{O}(f_{\alpha B}\varepsilon^{5})\nonumber\\
&&{}+{O}(n \gamma
\ln
B)+{O}(\varepsilon^{2} n \ln B )
\\
&=&n\sum_{\alpha\in S} \frac{f_{\alpha B}}{B}\psi_{\alpha}(\ln
n)+{O}(\varepsilon n).\nonumber
\end{eqnarray}
This proves the claim in the lattice case with $\varphi_B$ defined by
%
\begin{equation}\label{eqfixedcontinuous}
\varphi_B(q):=\sum_{\alpha\in S} \frac{f_{\alpha B}}B \psi_\alpha(q).
\end{equation}

It now only remains to prove that, although the functions $\psi_\alpha
( \cdot)$, $\alpha\in S$, are not continuous, the $d$-periodic
function $\varphi_B$ satisfies the bound in (\ref{eqcontinuityphiB}).
%
%
%
\begin{Lemma}
The function $\varphi_B$ defined in (\ref
{eqfixedcontinuous}) satisfies
\[
\sup_{|q-q'|\le\varepsilon^3} |\varphi_B(q)-\varphi_B(q')|\le K
\varepsilon
\ln(1/\varepsilon).
\]
\end{Lemma}
\begin{pf}
From the expresssion for $\psi_\alpha$ in (\ref{hopp2,1}), we have
\begin{eqnarray*}
\varphi_B(q)
&=& \frac{bd}\mu\sum_{\alpha\in S} \frac{f_{\alpha B}} B \sum_{k\dvtx
kd\le q} e^{kd-q} {\mathbf P}\bigl(q-kd+\ln V\in\bigl[\ln(\alpha-\gamma), \ln
\alpha\bigr)\bigr)\\
&=&\frac{bd}\mu\sum_{k\dvtx kd\le q} e^{kd-q} \sum_{\alpha\in S} \frac
{f_{\alpha B}} B {\mathbf P}\bigl(q-kd+\ln V\in\bigl[\ln(\alpha-\gamma), \ln
\alpha\bigr)\bigr).
\end{eqnarray*}
Note that, since $\gamma=\varepsilon^3$ and $\alpha\ge\varepsilon^2$,
\[
|{\ln}(\alpha-\gamma)-\ln\alpha|\sim\frac\gamma\alpha
\]
as $\varepsilon\to0$. As a consequence, for all $\varepsilon>0$ small enough,
the intervals involved in the definition of $\psi_\alpha$ satisfy,
uniformly in $\alpha\in S$,
\[
\frac{\varepsilon^3}2 <|{\ln}(\alpha-\gamma)-\ln\alpha|\le\varepsilon.
\]
In particular, since $\ln V\in d\mathbb Z$ almost surely, there is at
most one atom in the interval as soon as $\varepsilon<d$. It follows that,
if we choose $\delta=\varepsilon^3/2$, we have for any $q,q'$ such that
$|q-q'|< \delta$
\[
{\mathbf P}\bigl(q'-kd+\ln V\in\bigl[\ln(\alpha-\gamma), \ln\alpha\bigr)\bigr)
={\mathbf P}\bigl(q-kd+\ln V\in\bigl[\ln(\alpha'-\gamma), \ln\alpha'\bigr)\bigr)
\]
for some $\alpha'$ in $\{\alpha+\gamma,\alpha,\alpha-\gamma\}$. We
adopt the following point of view: for fixed $k$ and $q$, $S$ induces a
partition into the intervals $[q-kd-\ln(\alpha), q-kd-\ln(\alpha
-\gamma
))$, $\alpha\in S$. Each interval contains at most one atom of $-\ln
V$. Changing $q$ into $q'$ as above modifies the partition, but each
atom may only move to an adjacent interval. All atoms of $\ln V$ appear
in both sums, except if one is so far that it escapes the range of the
partition (recall that $\alpha\ge\varepsilon^2$). So following the atoms
of $-\ln V$ rather than the intervals in one or the other partition yields
\begin{eqnarray*}
\hspace*{-4pt}&&\frac{\mu}{bd}|\varphi_B(q)-\varphi_B(q')|\\
\hspace*{-4pt}&&\qquad\le\max_{x\in\{q,q'\}}\sum_{k\dvtx kd\le x} e^{kd - x +\delta} \sum
_{\alpha\in S} \max_{|\alpha'-\alpha|\le\gamma}\biggl| \frac
{f_{\alpha
B}} B -\frac{f_{\alpha' B}}B \biggr| \\
\hspace*{-4pt}&&\hspace*{184pt}{}\times{\mathbf P}\bigl(x-kd + \ln V \in
\bigl[\ln(\alpha-\gamma), \ln\alpha\bigr)\bigr)\\
\hspace*{-4pt}&&\qquad\quad{} + \max_{x\in\{q,q'\}} \sum_{k\dvtx kd\le x} e^{kd-x} \frac
{f_{\varepsilon^2B}} B,
\end{eqnarray*}
where the second term accounts for the escape of one atom. It follows that
\begin{eqnarray*}
\hspace*{-4pt}&&\frac{\mu}{bd}|\varphi_B(q)-\varphi_B(q')|\\
\hspace*{-4pt}&&\qquad\le\max_{x\in\{q,q'\}}\sum_{k\dvtx kd\le x} e^{kd - x +\delta} \sum
_{\alpha\in S} K \gamma\ln B \cdot{\mathbf P}\bigl(x-kd + \ln V \in\bigl[\ln
(\alpha-\gamma), \ln\alpha\bigr)\bigr)\\
\hspace*{-4pt}&&\qquad\quad{} + K \varepsilon^2 \ln B
\end{eqnarray*}
for some constant $K$, by Lemma~\ref{lemmatot} and the asymptotics for
$f_z$. Swapping the sums once again to recover the functions $\psi
_\alpha( \cdot)$, it follows that
\[
|\varphi_B(q)-\varphi_B(q')|
\le\frac{bd} \mu K \gamma e^\delta\ln B \cdot\sup_{x} \sum
_{\alpha
\in S} \psi_\alpha(x).
\]
However, since every summand is nonnegative, we have for any $x$
%
\begin{eqnarray}\label{eqboundsumalpha}\qquad
0&\le&\sum_{\alpha\in S} \psi_\alpha(x)
=\frac{bd}\mu\sum_{k\dvtx kd\le x} e^{kd-x} \sum_{\alpha\in S}
{\mathbf P}\bigl(x-kd +\ln
V\in\bigl[\ln(\alpha-\gamma),\ln\alpha\bigr)\bigr)\nonumber\\[-8pt]\\[-8pt]
&\le&\frac{bd} \mu\sum_{k\dvtx kd\le x} e^{kd-x}\le\frac{b
e^{d}}\mu.\nonumber
\end{eqnarray}
The desired bound follows: for any $q,q'$ such that $|q-q'|< \varepsilon
^3/2$, we have
\[
|\varphi_B(q)-\varphi_B(q')|\le K'' \varepsilon\ln(1/\varepsilon)
\]
for some constant $K''$ independent of $q,q'$ or
$\varepsilon$.
\end{pf}


%
%
%

\section{Extensions and concluding remarks}\label{secconcl}%

\subsection{An alternative notion of path length}

The notion of path length we have considered so far is the sum of the
depths of the \textit{items} in the tree. This is most natural when one
thinks about performance measures for algorithms or sorted data
structures. However, for some applications, it is sometimes important
to introduce a related notion of path length $\Upsilon(T)$, that is the
sum of the depths of \textit{nodes}:
\[
\Upsilon(T):=\sum_{u\in\mathcal U} |u| \one_{ \{ u\in T \} }=\sum
_{u\ne
\sigma} N_u,
\]
where $N_u$ denotes the number of nodes in the subtree rooted at $u$.
This notion of path length appears, for instance, in the analysis of
cutting-down processes. Suppose that you are given a rooted tree $T$.
Initially, the process starts with $T$. At each time step, a uniformly
random edge is cut, the portion of the tree that is disconnected from
the root is lost, and the process continues with the portion containing
the root. How many random cuts does it take to isolate the root? The
question originates in the seminal work of Meir and Moon
\cite{MeMo1970,MeMo1974}.\vadjust{\goodbreak}
Recently, the subject has regained interest, and new results have been
proved about the weak limit of the number of cuts when the initial tree
is randomly picked according to various distributions. See
\cite{Janson2006c,Holmgren2010b,Holmgren2010a,DrIkMoRo2009,IkMo2007} for
more references and details about the precise models and results.

For instance, Holmgren~\cite{Holmgren2010a} has proved that, when the initial
tree is a~split tree satisfying two general conditions (one on $\mathbf
{E}[\Upsilon(T^n)]$ and one on the number of nodes), the normalized number
of cuttings converges in distribution to a weakly 1-stable law (Theorem
1.1 there). Our Theorem~\ref{expectedpath} allows us to prove that one
of the conditions assumed in~\cite{Holmgren2010a} actually implies the
other. More precisely, the conditions assumed in~\cite{Holmgren2010a}
are that $\Upsilon(T^n)$ (the path length of nodes) satisfies
\[
\mathbf{E}[\Upsilon(T^n)]=\frac\alpha\mu n \ln n + \zeta n + o(n),
\]
and that the number of nodes $N=|T^n|$ verifies, for some constants
$\alpha>0$ and $\varepsilon>0$,
%
\begin{equation}\label{eqasymptoticsN}
\mathbf{E}[N]=\alpha n + f(n) \qquad\mbox{where } f(n)=O\biggl(\frac
n {\ln^{1+\varepsilon}n}\biggr).
\end{equation}
%
We deduce from Theorem~\ref{expectedpath}:
%
\begin{cor}\label{christmas}Suppose that $\ln V$ is nonlattice, and
assume that (\ref{eqasymptoticsN})
holds true; then, as $n\to\infty$,
\[
\mathbf{E}[\Upsilon(T^n)]=\frac\alpha\mu n \ln n + \zeta n + o(n).
\]
\end{cor}
\begin{Remarks*} The assumption in (\ref{eqasymptoticsN})
is just slightly stronger than the estimate proved by Holmgren
\cite{Holmgren2009}, that is, that for split tree with nonlattice
$\ln
V$, we have $f(n)=o(n)$. Moreover, the assumption in
(\ref{eqasymptoticsN}) does make sense, since it is known to hold,
for instance, for $m$-ary search trees
\cite{Knuth1973b,Baeza1987,MaPi1989,ChPo2004}: for such random trees, $f(n)$
is $o(\sqrt{n})$ when $m\leq26$ and is ${O}(n^{1-\varepsilon})$ when
$m\geq27$. On the other hand, it is also known that the condition
in~(\ref{eqasymptoticsN}) does not always hold. For instance,
Flajolet et al.~\cite{FlRoVa2010a} proved that, in the case of binary
tries generated
by a memoryless source with probabilities $p_1, p_2$ such that $(\log
p_1)/(\log p_2)$ is a Liouville number, then the error term $f(n)$ can
come arbitrarily close to $O(n)$ [but of course, stays $o(n)$].
See~\cite{FlRoVa2010a}, page~249, and the monograph by Baker~\cite{Baker1990a}
for more information about Liouville numbers.
\end{Remarks*}
\begin{pf*}{Sketch of proof}
Define $q(n)$ and $r(n)$ by
\[
\mathbf{E}[\Psi(T^n)]=\frac1 \mu n \ln n + n q(n) \quad\mbox{and}\quad
\mathbf{E}[\Upsilon(T^n)]=\frac\alpha\mu n \ln n + n r(n).
\]
Let $\Delta_n:=\alpha nq(n)-nr(n)$, and note that
%
\begin{equation}\label{snow}
\Delta_n=\alpha\mathbf{E}[\Psi{(T^{n})}]-\mathbf{E}[\Upsilon{(T^{n})}].\vadjust{\goodbreak}
\end{equation}
Since, by Theorem~\ref{expectedpath}, $q(n)$ converges as $n\to\infty$,
it suffices to prove that $\Delta_n/n$ also converges to some constant.
From (\ref{snow}) and the assumption in~(\ref{eqasymptoticsN}) we obtain
%
\begin{eqnarray}\label{snowy}
\Delta_n&=&\alpha\mathbf{E}\biggl[\sum_{v\neq\sigma} n_v\biggr]-\mathbf
{E}\biggl[\sum_{v\neq\sigma} \biggl(\alpha n_v+{O}\biggl(\frac{n_v}{\ln
^{1+\varepsilon} n_v}\biggr)\biggr)\biggr]\nonumber\\[-8pt]\\[-8pt]
&=&\mathbf{E}\biggl[\sum_{v} {O}\biggl(\frac
{n_v}{\log^{1+\varepsilon}n_v}\biggr)\biggr].\nonumber
\end{eqnarray}
[The constants hidden in the $ O(\cdot)$ above are the same for every
term.]

Consider the subtrees $T_{r}$, $r\in R$, introduced in the course of
the proof of Theorem~\ref{expectedpath}. Recall that a node $r$ is in
$R$ if it is the first on its path from the root such that $n L_r \le
B$, for some parameter $B$. In the following, we take $B=\delta^{-8}$,
for $\delta>0$. We now show that the main contribution to~$\Delta_n$ is
accounted for by the nodes in the subtrees $T_r$, $r\in R$; in other
words $\Delta_n=\mathbf{E}[\sum_{r\in R} \Delta_{n_r}]+o(n)$, where
\[
\Delta_{n_r}
=\alpha\mathbf{E}[\Psi{(T_{r})}|n_r]-\mathbf{E}[\Upsilon{(T_{r})}|n_r].
\]
To see this, observe that we deduce from (\ref{snowy}) and
(\ref{eqasymptoticsN}) that
\begin{eqnarray*}
\Delta_n - \mathbf{E}\biggl[\sum_{r\in R} \Delta_{n_r}\biggr]
&=&\mathbf{E}\biggl[\mathop{\sum_{v\notin T_{r}, r\in R,}}_{{v\neq\sigma
}}{O} \biggl(\frac{n_v}{\log^{1+\varepsilon}n_v}\biggr)\biggr]\\
&=&\mathbf{E}\biggl[\sum_{k\ge0}\mathop{\sum_{v\notin T_{r}, r\in
R,}}_{{2^{k}\leq n_v<2^{k+1}}} {O}\biggl(\frac{n_v}{\log^{1+\varepsilon
}n_v}\biggr)\biggr]+ O\biggl(\frac n {\log n}\biggr).
\end{eqnarray*}
We split the sum in $k$ above at some constant $K$ to be chosen later.
By Lem\-ma~\ref{cardinality} and since the expected number of nodes
$v\in
T^n$ with $nL_v\ge B$ is $ O(n/B)$, we obtain
\begin{eqnarray*}
\Delta_n - \mathbf{E}\biggl[\sum_{r\in R} \Delta_{n_r}\biggr]
&=&\sum_{k> K} {O}\biggl(\frac{n}{2^{k}}\cdot\frac{ 2^{k}}{
k^{1+\varepsilon
}}\biggr) + \sum_{0\le k \le K} {O}\biggl(\frac{n}{B}\cdot\frac{ 2^{k}}{
k^{1+\varepsilon}}\biggr)+o(n)\\
&=& O(n K^{-\varepsilon})+ O(n K 2^K/B)+o(n).
\end{eqnarray*}
We choose $K=\lfloor a\ln(1/\delta) \rfloor$, for some small
constant $a>0$.
Since $\delta>0$ was arbitrary, the claim follows.

Now since $\Delta_{n_r}={O}(n_r\ln n_r)$, the proof of
Proposition \ref
{lemma2} (in the nonlattice case) may be extended to show that
$\mathbf{E}[\sum_{r\in R}\Delta_{n_r}]=n\zeta+o(n)$ for some
constant~$\zeta$. The
details are omitted.
\end{pf*}

\subsection{Beyond split trees and multinomial partitions}

To conclude, we indicate the lines of the arguments to extend the
applicability of our main theorem to a greater family of random trees.\vadjust{\goodbreak}
The model of split trees~\cite{Devroye1998} supposes that the
distribution of the subtree cardinalities $n_1,n_2,\ldots, n_b$ of
a~node of cardinality $n$ is \textit{exactly} of the form
%
\begin{equation}\label{eqmultform}\quad
(n_1,n_2,\ldots, n_b)=\operatorname{Mult}(n-s_0-bs_1,V_1,V_2,\ldots,
V_b)+(s_1,s_1,\ldots, s_1)
\end{equation}
for a random vector $(V_1,\ldots, V_b)$; in particular, the vector
$(V_1,\ldots, V_b)$ cannot depend on $n$. Although many important data
structures satisfy this property, some other more combinatorial
examples do not; see, for instance, the case of increasing trees
\cite{BeFlSa1992}.

Also, the reader might have noticed that our proof does not quite use
the full strength of the assumption in (\ref{eqmultform}). Indeed,
our proof mainly uses two facts: first, that the sequence of subtree
sizes along a branch is well approximated by the product form $nL_u =n
\prod_{v\preceq u} V_v$, which modulo some details about $C(\mathcal
V)$, implies that
%
\[
X\eqdist\sum_{k=1}^b V_k X^{(k)} + C(\mathcal V);
\]
and second, that the addition of some items to the tree only modifies
moderately $\mathbf{E}[\Psi(T)]$ (see Lemma~\ref{lemmatot}).

The two requirements are satisfied when the items are distributed in
subtrees according to (\ref{eqmultform}).
We now indicate why our result would still hold under the much weaker
condition that there exists a vector $\mathcal V=(V_1,\ldots, V_b)$ such
that the cardinalities $n_1,\ldots, n_b$ of the children of a node of
cardinality $n$ satisfy
%
\begin{equation}\label{eqlimitv}
\biggl(\frac{n_1}n, \frac{n_2} n,\ldots, \frac{n_b}n\biggr)\to
(V_1,V_2,\ldots
, V_b) \qquad\mbox{in distribution}
\end{equation}
as $n\to\infty$. Of course, the copies of the limit vectors $\mathcal
V$ at distinct nodes should be independent. The general shape of trees
under this model has recently been completed by work by Broutin
et al.~\cite{BrDeMcSa2008} (see also
Drmota~\cite{Drmota2009a} who treats the model of increasing trees by
Bergeron et al.~\cite{BeFlSa1992} more directly). 

One should be easily convinced that the relaxed condition in (\ref
{eqlimitv}) should be sufficient for the result to hold:
\begin{itemize}
\item Proposition~\ref{lemma1} may be extended using the coupling
arguments already used in~\cite{BrDeMcSa2008}, proving that the
contribution of the top of the tree to the path length may be estimated
using renewal functions associated to the limit vector $\mathcal V$.
\item Similarly, the extension of Proposition~\ref{lemma2} relies on
the same coupling argument (the overshoot there is still approximated
by that of the limit vector). Here, it is important to note that the
proof of smoothness of the path length (Lemma~\ref{lemmatot}) requires
the existence of a fixed function $g$ such that the size $|T^n|$ of a
``generalized'' split tree of cardinality $n$ satisfies $|T^n|\le g(n)$
with probability 1 (at least our proof does). This was already
necessary for the results on the shape of the trees in
\cite{BrDeMcSa2008} to hold. The constraint is not too strong, since it
holds as soon as $s_0$ or $s_1$ is nonzero,\vadjust{\goodbreak} and any function would do,
regardless of its growth. (This is another reason why the case of
digital trees should be treated separately: for such trees, the size of
a tree containing two items can be arbitrarily large.)
\item As already noted in Section~\ref{seccontraction}, the part of
the proof relative to the contraction method in
\cite{NeRu1999,NeRu2004} will go through as long as the coefficients
$C_n(\overline n)$ converge, and the expansion for mean implies their
convergence.
\end{itemize}

\section*{Acknowledgments}

We would like to warmly thank Svante Janson and Ralph Neininger for
very helpful discussions. We are also grateful to an anonymous referee
for his helpful comments.


%

\printaddresses

\end{document}